\documentclass[11pt]{article}
\usepackage[T1]{fontenc}
\usepackage{amsmath,amsthm}         % load before newtxmath (it defines \openbox etc.)
\usepackage{newtxtext}
\usepackage[varvw]{newtxmath}       % selects Times Roman as basic font
\usepackage{thmtools}
\usepackage{thm-restate}
\usepackage{graphicx}
\usepackage{enumerate}
\usepackage{bm}
\usepackage{float}
\usepackage{cite}
\usepackage[margin=1in]{geometry}
\usepackage{hyperref}
\usepackage{color}
\hypersetup{
    colorlinks=true,
    linkcolor=black,
    filecolor=magenta,
    urlcolor=black,
    citecolor=black,
}
\theoremstyle{plain}
\newtheorem{theorem}{Theorem}[section]
\newtheorem{corollary}[theorem]{Corollary}

\newtheorem{lemma}[theorem]{Lemma}

\theoremstyle{definition}
\newtheorem{definition}[theorem]{Definition}

\theoremstyle{remark}
\newtheorem{remark}[theorem]{Remark}
\newtheorem{example}[theorem]{Example}
\newtheorem*{theorem*}{Theorem}
\newcommand{\Z}{\ensuremath{\mathbb{Z}}}
\newcommand{\Q}{\ensuremath{\mathbb{Q}}}
\newcommand{\R}{\ensuremath{\mathbb{R}}}
\newcommand{\N}{\ensuremath{\mathbb{N}}}
\newcommand{\T}{\ensuremath{\mathsf{T}}}
\makeatletter
\renewenvironment{proof}[1][]{\par
  \pushQED{\qed}%
  \normalfont \topsep6\p@\@plus6\p@\relax
  \trivlist
  \item[\hskip\labelsep\itshape
    \proofname\if\relax\detokenize{#1}\relax\else\ (#1)\fi\@addpunct{.}]\ignorespaces
}{%
  \popQED\endtrivlist\@endpefalse
}
\makeatother
\newcommand{\keywords}[1]{%
  \par\addvspace{\baselineskip}\noindent
  \textbf{Keywords:} #1\par}

\begin{document}
\title{Frank--Wolfe Theorems for Mixed Integer Quadratic Optimization}
\author{Alberto Del Pia\thanks{%
  Department of Industrial and Systems Engineering \& Wisconsin Institute for Discovery,
  University of Wisconsin--Madison, Madison, WI, USA.
  \texttt{delpia@wisc.edu}}
  \and
  Brandon Sommerfeld\thanks{%
  Department of Computer Sciences, University of Wisconsin--Madison, Madison, WI, USA.
  \texttt{bdsommerfeld@wisc.edu}}%
}

\date{November 10, 2026}

\maketitle              % typeset the header of the contribution

\begin{abstract}
The Frank--Wolfe theorem states that a quadratic function bounded from below
on a polyhedron attains its infimum. We study when this remains true in the presence of
integer variables, for sets defined by quadratic inequalities. Our main result is a complete
classification for sets defined by convex quadratic inequalities with rational coefficients
and a rational affine hull, a condition that holds for all full-dimensional sets. One integer
variable never destroys attainment: if every rational quadratic objective that is bounded from
below attains its infimum over the set, even after adding rational linear constraints, 
the same holds with one integer variable. In dimension
at most 2, integer variables never destroy attainment. In dimension at least 3 and with two or
more integer variables, attainment holds for every such objective if and only if the set has
bounded curvature, that is, it is bounded except along directions in which all constraints are
linear. Without convexity, these results remain valid, except that bounded curvature is no
longer necessary. Moreover, integer variables can destroy attainment only if there are at least
two of them and at least one continuous variable. The result for one integer variable follows
from a lifting theorem, which also extends known results for cubic objectives and convex
polynomial constraints. Finally, for convex constraints and rational data, attainment holds when the objective is
bounded from above, convex, or quasi-convex over a full-dimensional convex set containing the
feasible region.

\end{abstract}

\section{Introduction}

A classical theorem of Frank and Wolfe~\cite{FrankWolfe} states that a quadratic function that is
bounded from below on a polyhedron attains its infimum. Results of this kind are called
\emph{Frank--Wolfe type theorems}, and many of them are known for continuous problems, see
Section~\ref{sec:related}. In this paper, we ask when they remain true in the presence of integer variables.

\paragraph{Setting.}
We consider mixed integer quadratically constrained quadratic programs of the form
\begin{align}
\label{prob:MIQCQP}
\begin{split}
    \inf\; &f_{0}(x) \\
    \text{s.t.}\; & x \in S \cap P \cap \Lambda.
\end{split}
\end{align}
Here
\begin{align}
\label{def:feasible_set}
    S := \{x \in \R^{n} : f_{i}(x) \leq 0, \quad i=1,2,\dots,m\},
\end{align}
where, for each $j = 0,1,\dots,m$,
\[
    f_{j}(x) = x^{\T} Q_{j} x + q_{j}^{\T} x + c_{j},
\]
with a symmetric matrix $Q_{j} \in \R^{n \times n}$, $q_{j} \in \R^{n}$, and $c_{j} \in \R$.
Moreover, $P \subseteq \R^{n}$ is a polyhedron, and $\Lambda$ is a \emph{mixed integer lattice}, that is,
a set of the form
\[
    \Lambda = \tau\left(\Z^{p} \times \R^{n-p}\right),
\]
where $\tau(x) = Mx + \bar{x}$ for an invertible $M \in \R^{n \times n}$ and $\bar{x} \in \R^{n}$.
We say that $\Lambda$ has \emph{$p$ integer components}. For $p = 0$, we obtain the continuous case
$\Lambda = \R^{n}$. We say that~\eqref{prob:MIQCQP} \emph{attains its infimum} if there exists
$x \in S \cap P \cap \Lambda$ such that $f_{0}(x) \leq f_{0}(y)$ for every $y \in S \cap P \cap \Lambda$.
In particular, the infimum is then finite. To avoid trivial cases, we assume throughout the paper that
$S \cap P \cap \Lambda$ is nonempty.

An instance of~\eqref{prob:MIQCQP} has \emph{rational data} if all $f_{j}$ have \emph{rational coefficients},
that is, all entries of $Q_{j}$, $q_{j}$, and $c_{j}$ are rational, $P$ is \emph{rational}, that is, defined by
finitely many linear inequalities with rational coefficients, and $\Lambda$ is \emph{rational}, that is,
$M$ and $\bar{x}$ can be chosen rational. In this introduction, we always assume rational data.

\paragraph{The Frank--Wolfe property.}
Whether an optimization problem has an optimal solution is one of the most basic questions in
optimization. Algorithms, optimality conditions, and complexity results often assume that the infimum
is attained. For mixed integer quadratic optimization, this leads to a fundamental question: for which
sets $S$ is the infimum in~\eqref{prob:MIQCQP} attained for every $f_{0}$, $P$, and $\Lambda$ for which
it is finite? For continuous variables, that is, if only $\Lambda = \R^{n}$ is allowed,
the Frank--Wolfe theorem shows that $S = \R^{n}$ has this property, and Frank--Wolfe type theorems give
other such sets. This question leads to the following definition.
\begin{definition}
    A set $S \subseteq \R^{n}$ has the \emph{Frank--Wolfe property $\text{FW}_{p}$}, or simply is
    \emph{$\text{FW}_{p}$}, if, for every rational mixed integer lattice $\Lambda$
    with at most $p$ integer components, every rational polyhedron $P$, and every quadratic function
    $f_{0}$ with rational coefficients that is bounded from below on $S \cap P \cap \Lambda \neq \emptyset$,
    problem~\eqref{prob:MIQCQP} attains its infimum.
\end{definition}
The intersection with $P$ makes the property hereditary: if $S$ is $\text{FW}_{p}$,
then so is $S \cap P$ for every rational polyhedron $P$. Moreover, $\text{FW}_{p}$ implies
$\text{FW}_{p-1}$, and $\text{FW}_{0}$ only concerns continuous variables. In this language,
the Frank--Wolfe theorem shows that $\R^{n}$ is $\text{FW}_{0}$, and hence so is every rational polyhedron.

\subsection{Main result}
\label{sec:main_result}

The following notion turns out to characterize the Frank--Wolfe property for convex constraints.
\begin{definition}
    Let $S$ be defined as in~\eqref{def:feasible_set}, and let $N := \bigcap_{i=1}^{m} \ker Q_{i}$.
    We say that $S$ has \emph{bounded curvature} if $S \subseteq B + N$ for a bounded set $B$.
\end{definition}
In words, $S$ is bounded except along the directions in $N$, along which all constraints are linear.
Equivalently, the orthogonal projection of $S$ onto $N^{\perp} = \sum_{i=1}^{m} \operatorname{range} Q_{i}$ is bounded.
For example, polyhedra and compact sets have bounded curvature, while the paraboloid
$\{x \in \R^{3} : x_{1}^{2} + x_{2}^{2} \leq x_{3}\}$ does not.

Our main result is the following complete classification.
\begin{theorem} \label{thm:main}
    Let $S$ be given by~\eqref{def:feasible_set}, where all $f_{i}$ are convex with rational coefficients,
    and assume that the affine hull of $S$ is rational. Then
    \begin{enumerate}[(i)]
        \item $S$ is $\text{FW}_{1}$ if and only if $S$ is $\text{FW}_{0}$.
        \item If the affine hull of $S$ has dimension at most 2, then $S$ is $\text{FW}_{n}$ if and only if
            $S$ is $\text{FW}_{0}$.
        \item If the affine hull of $S$ has dimension at least 3, then for every $2 \leq p \leq n$,
            the set $S$ is $\text{FW}_{p}$ if and only if $S$ has bounded curvature.
    \end{enumerate}
\end{theorem}
In words, one integer variable never destroys attainment, and in dimension at most 2, no number of
integer variables does. This follows from Theorem~\ref{thm:generic_one_integer}, 
which lifts attainability results from the purely continuous setting to the one integer setting.
In dimension at least 3 and with two or more integer variables, attainment for all
quadratic objectives is equivalent to bounded curvature. The assumption on the affine hull
holds for full-dimensional sets, and for sets defined by one convex quadratic inequality and linear
inequalities~\cite[Theorem 1.2]{DelPia24}. 

The threshold 3 in (iii) is sharp. The set $\{x \in \R^{2} : x_{1}^{2} \leq x_{2}\}$ is $\text{FW}_{0}$
by~\cite[Theorem 2]{LuoZhang}, and hence $\text{FW}_{2}$ by (ii), but it does not have bounded curvature.
For a full-dimensional set in $\R^{3}$ as in Theorem~\ref{thm:main}, the theorem leaves three
possibilities: the set is
$\text{FW}_{3}$, it is $\text{FW}_{1}$ but not $\text{FW}_{2}$, or it is not $\text{FW}_{0}$.
All three occur. Sets with bounded curvature are $\text{FW}_{3}$. The set $\{x \in \R^{3} : x_{1}^{2} + x_{2}^{2} \leq x_{3}\}$ is $\text{FW}_{1}$
by~\cite[Theorem 2]{LuoZhang} and (i), but it is not $\text{FW}_{2}$ by Example~\ref{ex:multiple_integer_convex}.
Finally, the set $\{x \in \R^{3} : x_{2}^{2} \leq x_{1}, \ x_{3}^{2} \leq x_{1} + 1\}$ is not even
$\text{FW}_{0}$~\cite[Section 3]{BelousovKlatte02}.

\subsection{Extensions and consequences}
\label{sec:further}

We now discuss which parts of Theorem~\ref{thm:main} extend to nonconvex constraints, what our results
imply about the role of integer variables, how conditions on the objective can replace bounded curvature,
and what happens without rational data. We conclude by relating our results to those of Luo and Zhang~\cite{LuoZhang}.

\paragraph{Nonconvex constraints.}
Without convexity, parts (i) and (ii) of Theorem~\ref{thm:main} remain valid, see
Corollaries~\ref{cor:quadratic_one_integer} and~\ref{cor:low_dimension}. In part (iii), bounded curvature
remains sufficient for $S$ to be $\text{FW}_{n}$, see Corollary~\ref{cor:forward}, but it is no longer
necessary. Remark~\ref{rmk:nonconvex} gives a nonconvex set that is $\text{FW}_{n}$ but does not have
bounded curvature.

\paragraph{The role of integer variables.}
Let $S$ be given by~\eqref{def:feasible_set}, convex or not, and assume that $S$ is $\text{FW}_{0}$,
so that attainment holds with continuous variables only. Then integer variables can destroy attainment
only if there are at least two of them, by Corollary~\ref{cor:quadratic_one_integer}, and at least one
continuous variable, since attainment always holds when all variables are integer, see
Corollary~\ref{cor:all_integer}. Moreover, this can only happen if $S$ does not have bounded curvature,
by Corollary~\ref{cor:forward}.

\paragraph{Conditions on the objective.}
Theorem~\ref{thm:main} concerns only the set $S$, and it quantifies over all
mixed integer lattices and all quadratic objectives. Under additional assumptions on
the objective, the infimum is attained even when $S$ does not have bounded curvature, as the following
result shows. 
\begin{theorem} \label{thm:bounded_above}
    Consider an instance of~\eqref{prob:MIQCQP} with rational data, where all $f_{i}$ are convex
    and $f_{0}$ is bounded from above and from below on $S \cap P \cap \Lambda$.
    Then $f_{0}$ attains its infimum on $S \cap P \cap \Lambda$.
\end{theorem}
We prove Theorem~\ref{thm:bounded_above} as Theorem~\ref{thm:convex} in Section~\ref{sec:convex}.
Together with its consequences for convex and quasi-convex objectives, Theorem~\ref{thm:bounded_above}
shows that, for convex constraints, the objective in any nonattainment
example is nonconvex and unbounded from above, see Remark~\ref{rmk:tightness}.

\paragraph{Irrational data.}
So far we assumed rational data. With one integer variable, the lifting theorem in
Section~\ref{sec:one_integer} does not need rational data. With two integer variables, rational data
are needed, even for polyhedra, see Appendix~\ref{sec:irrational}, where Table~\ref{tab:summary}
summarizes both settings.

\paragraph{Relation to Luo and Zhang.}
Our results extend to the mixed integer setting the results of Luo and Zhang~\cite{LuoZhang}
on attainment in continuous quadratically constrained quadratic programming.
First, they showed that a set $S$ defined by a single convex quadratic inequality and
linear inequalities is $\text{FW}_{0}$~\cite[Theorem 2]{LuoZhang}. By Corollary~\ref{cor:quadratic_one_integer},
such a set is always $\text{FW}_{1}$, and if $S$ has dimension at most 2 and a rational affine hull, then 
it will be $\text{FW}_{n}$ by Corollary~\ref{cor:low_dimension}. For $S$ 
which has dimension at least 3 and a rational affine hull, $S$ is $\text{FW}_{p}$ 
for any $2 \leq p \leq n$ if and only if $S$ has bounded curvature by Corollary~\ref{cor:characterization}.
Additionally, for a single convex quadratic inequality $x^{\T} Q x + q^{\T} x + c \leq 0$ with
$Q \neq 0$ and $n \geq 3$, Corollary~\ref{cor:one_constraint} shows that $S$ is
$\text{FW}_{n}$ if and only if $q \in \operatorname{range} Q$.

Second, Luo and Zhang showed that the infimum is attained when all constraints are convex
and $f_{0}$ is quasi-convex over the polyhedron defined by the linear
constraints~\cite[Theorem 3]{LuoZhang}. Corollary~\ref{cor:quasi-convex} shows that, for any number of
integer variables, it suffices that $f_{0}$ is convex or quasi-convex over
some full-dimensional convex set containing $S$. For at most one integer variable,
Corollary~\ref{cor:one_integer_quasi-convex} only requires that $f_{0}$ is quasi-convex over $S$.
Luo and Zhang also showed that nonconvex constraints can
define sets that are not $\text{FW}_{0}$, even when the feasible region is
convex~\cite[Example 3]{LuoZhang}.

\subsection{Related literature}
\label{sec:related}

Frank--Wolfe type theorems for continuous problems have been studied extensively. Alternative proofs
and extensions of the Frank--Wolfe theorem are given in~\cite{Eaves71,BlumOet72,Perold80}. Attainment is
also known for cubic objectives with linear constraints~\cite{AndBelShi82,Klatte19}, for convex polynomial
objectives and constraints~\cite{Belousov77,BelousovKlatte02}, and for other classes of objectives and
constraints~\cite{Obuchowska06,NghiTam25,DinhHaPham13}. Ozdaglar and Tseng~\cite{OzdaglarTseng06} and
Bertsekas and Tseng~\cite{BertsekasTseng07} gave unified existence results based on recession and
asymptotic directions, and Bertsekas and Tseng mention extensions with integer variables as a topic for
future research~\cite[Section 4]{BertsekasTseng07}.

Other works characterize the feasible sets on which Frank--Wolfe type theorems hold.
Kummer~\cite{Kummer77} proved such a theorem for sets of the form $K + D$, where $K$ is compact and $D$
is a polyhedral cone. Martínez-Legaz et al.~\cite[Theorem 12.2]{Martinez19} proved the converse: if every
quadratic function that is bounded from below on $K + D$ attains its infimum there, for a closed convex
cone $D$, then $D$ is polyhedral. This is a continuous analog of Theorem~\ref{thm:counterexample}. They
also characterize the convex sets on which every quadratic objective that is quasi-convex on the set and
bounded from below attains its infimum~\cite[Theorem 3.1]{Martinez18}. Corollary~\ref{cor:quasi-convex}
is a result of this type in the mixed integer setting.

In the mixed integer setting, Meyer~\cite{Meyer74} proved attainment for linear objectives and linear
constraints with rational coefficients, and Bank and Hansel~\cite{BanHan84} for quadratic objectives with
linear constraints. Bank and Mandel~\cite{BankMandel88} extended this result to
sets of the form $K + V$, where $K$ is compact and $V$ is generated by finitely many integer vectors.
Corollary~\ref{cor:forward} follows from their result, and the converse in Theorem~\ref{thm:main}(iii)
is new. For integer programs with convex or reverse convex constraints,
Obuchowska~\cite{Obuchowska08,Obuchowska10} studied boundedness and the existence of optimal solutions.

\paragraph{Organization of the paper.}
Section~\ref{sec:notation} introduces our notation. Section~\ref{sec:one_integer} lifts attainment results
from continuous variables to one integer variable, and proves part~(i) of Theorem~\ref{thm:main}.
Section~\ref{sec:bounded_curvature} considers two or more integer variables and proves the rest of
Theorem~\ref{thm:main}.
Section~\ref{sec:convex} considers conditions on the objective and proves Theorem~\ref{thm:bounded_above}.
The appendix contains the omitted proofs,
further consequences of the lifting theorem, and the results on irrational data.

\section{Preliminaries}
\label{sec:notation}

Throughout the paper, $f_{j}$, $Q_{j}$, $q_{j}$, and $c_{j}$ always refer to an instance
of~\eqref{prob:MIQCQP}. Results that are more general than~\eqref{prob:MIQCQP} use other notation,
typically $g_{j}$. Since a linear inequality is a quadratic inequality with $Q_{i} = 0$, in results about
a single instance of~\eqref{prob:MIQCQP} we include the inequalities defining $P$ among those defining $S$,
and write $S \cap \Lambda$ for the feasible region.

We say that a direction $d \in \R^{n}$ is a \emph{recession direction} of an inequality $f(x) \leq 0$ 
where $f : \R^{n} \to \R$,
if $f(x+td) \leq 0$ for every $x$ with $f(x) \leq 0$ and $t \geq 0$.
Similarly, $d$ is a \emph{recession direction} of a set $X \subseteq \R^{n}$ if 
\[\forall x \in X, \forall t \geq 0, x+td \in X .\]
We denote by $\operatorname{rec}(X)$ the set of recession directions of $X$.
The \emph{dimension} of a set is the dimension of its affine hull.
We denote by $e^{i}$ the $i$-th standard unit vector, whose dimension is clear from the context.
An affine subspace of $\R^{n}$ is \emph{rational} if it can be defined by linear equations
with rational coefficients, or equivalently, if it is of the form $\{x_{0} + Wz : z \in \R^{k}\}$
with $x_{0} \in \Q^{n}$ and $W \in \Q^{n \times k}$.

Throughout the paper, when we refer to a quadratic function $g : \R^{n} \to \R$, the canonical 
form is $g(x) = x^{\T} Q x + q^{\T} x + c$, where $Q \in \R^{n \times n}$ is a symmetric matrix, 
$q \in \R^{n}$, and $c \in \R$. The function $g$ is convex if and only if $Q$ is positive semidefinite~\cite[Theorem 4.5]{Rockafellar97}.

\begin{remark} \label{rmk:extension}
    All examples in this paper can be extended by any number of additional integer
    or continuous variables that appear neither in the constraints nor in the objective.
    Hence they generalize to higher dimensions.
\end{remark}

Our proofs use four elementary lemmas on recession directions and rational approximation.
We state and prove them in Appendix~\ref{app:prelim}.

\section{One integer variable} 
\label{sec:one_integer}

In this section, we show that attainability in the continuous 
setting can be lifted to the single integer setting. 
Specifically, we consider problems of the form  

\begin{align}
\label{prob:generic_one_integer}
\begin{split}
    \inf\; & g_{0}(x) \\
    \text{s.t.}\; & g_{i}(x) \leq 0, \quad i=1,2,\dots,m, \\ 
    & x \in \Lambda,
\end{split}
\end{align}
where for $j=0,1,\dots, m$, $g_{j}: \R^{n} \to \R$ is a general function, 
and $\Lambda = M\left(\Z \times \R^{n-1}\right) + \bar{x}$
for an invertible $M \in \R^{n \times n}$ and $\bar{x} \in \R^{n}$. Thus $\Lambda$ is a mixed integer lattice
with one integer component, and in this section it need not be rational. Let $a^{*} := M^{-\T} e^{1}$, so that
$\Lambda = \{x \in \R^{n} : (a^{*})^{\T}(x - \bar{x}) \in \Z\}$.

We also utilize the related continuous problem with additional constraints, given by
\begin{align}
\label{prob:relaxed_with_linear}
\begin{split}
    \inf\; & g_{0}(x) \\
    \text{s.t.}\; & g_{i}(x) \leq 0, \quad i=1,2,\dots,m, \\ 
    &L \leq a^{\T} x \leq U \\
    &x \in \R^{n},
\end{split}
\end{align}
obtained by relaxing the mixed integer lattice in~\eqref{prob:generic_one_integer}
and adding bounds $L,U \in \R \cup \{-\infty, \infty\}$ on $a^{\T} x$
for some $a \in \R^{n}$, but keeping the same functions $g_{j}$ for $j=0,1,\dots,m$.

Attaining the infimum is defined as for~\eqref{prob:MIQCQP}.

\begin{theorem} \label{thm:generic_one_integer}
    Consider a problem of the form~\eqref{prob:generic_one_integer} with finite infimum $c$,
    and assume that the set
    \[\{x \in \R^{n} : g_{0}(x) \leq c, g_{i}(x) \leq 0, \quad i=1,2,\dots, m \}\]
    has finitely many connected components.
    Assume that for $a = a^{*}$ and every choice of $L,U \in \left((a^{*})^{\T}\bar{x} + \Z\right) \cup \{-\infty, \infty\}$
    for which~\eqref{prob:relaxed_with_linear} has a finite infimum, this infimum is attained.
    Then~\eqref{prob:generic_one_integer} attains its infimum.
\end{theorem}

The idea of the proof is as follows. Consider a sequence of feasible points whose objective values
converge to the infimum. If their integer coordinates are bounded, then only finitely many slices with
a fixed integer coordinate are involved, and each slice attains its infimum by assumption. Otherwise,
the sublevel set of the continuous relaxation at the infimum has finitely many connected components,
and each of them lies strictly between two consecutive values of the integer coordinate. Beyond the
last of them, a relaxation with one additional linear constraint has the same infimum but does not
attain it.

We prove Theorem~\ref{thm:generic_one_integer} in Appendix~\ref{app:one_integer}.

We now prove part~(i) of Theorem~\ref{thm:main}, for every set given by~\eqref{def:feasible_set}.

\begin{corollary} \label{cor:quadratic_one_integer}
    For any $S$ defined as in~\eqref{def:feasible_set},
    $S$ is $\text{FW}_{0}$ if and only if $S$ is $\text{FW}_{1}$. 
\end{corollary}

\begin{proof}
    By definition, every $\text{FW}_{1}$ set is $\text{FW}_{0}$, so we only need to prove the converse.

    Let $P$ be a rational polyhedron.
    Consider the set  
    \[S_{c} := \{x \in \R^{n}: f_{0}(x) \leq c, f_{i}(x) \leq 0, \quad i = 1,2,\dots,m \}. \]
    Milnor~\cite[Theorem 3]{Milnor} shows that any set defined by polynomial inequalities 
    has a finite sum of Betti numbers, and this also 
    applies to $S_{c} \cap P$ since $P$ is given by linear inequalities. 
    Since the 0th Betti number is the number of connected components, 
    this shows that for any $c$ the number of connected components of $S_{c} \cap P$ is finite.

    Since $S$ is $\text{FW}_{0}$, so is $S \cap P$. Mixed integer lattices with no integer components are
    covered by the definition of $\text{FW}_{0}$. Consider a rational mixed integer lattice $\Lambda$ with one
    integer component, and an instance of~\eqref{prob:generic_one_integer} with
    $\{x \in \R^{n} : g_{i}(x) \leq 0, \ i=1,2,\dots,m\} = S \cap P$
    and an objective $g_{0}$ that is quadratic with rational coefficients.
    Since $\Lambda$ is rational, $a^{*}$ and $(a^{*})^{\T}\bar{x}$ are rational. Hence each
    relaxation~\eqref{prob:relaxed_with_linear} in Theorem~\ref{thm:generic_one_integer} only adds rational
    linear inequalities, and it attains its infimum when this is finite. By Theorem~\ref{thm:generic_one_integer},
    problem~\eqref{prob:generic_one_integer} attains its infimum. Hence $S$ is $\text{FW}_{1}$.
\end{proof}

In the same way, Theorem~\ref{thm:generic_one_integer} lifts other attainment results from continuous
variables to one integer variable. In Appendix~\ref{app:one_integer}, we extend the results of Luo and
Zhang~\cite[Theorem 2]{LuoZhang} for one convex quadratic constraint, of Klatte~\cite{Klatte19} for cubic
objectives, and of Belousov and Klatte~\cite[Theorem 3]{BelousovKlatte02} for convex polynomials.
These extensions do not require rational data or a rational mixed integer lattice. In contrast,
rational data are needed with two integer variables, see Appendix~\ref{sec:irrational}.

\section{Two or more integer variables}
\label{sec:bounded_curvature}

In this section, we consider mixed integer lattices with at least two integer components.
We first show that bounded curvature implies the property $\text{FW}_{n}$, for every set
with rational coefficients. We then show that bounded curvature is necessary for convex sets of
dimension at least 3 with a rational affine hull. Finally, we treat lattices with $n$ integer
components and sets of dimension at most 2.

\begin{lemma} \label{lem:bounded_curvature}
    Let $S$ be defined as in~\eqref{def:feasible_set}, and assume that $S$ is nonempty.
    Let $N := \bigcap_{i=1}^{m} \ker Q_{i}$ and $R := \{w \in N : q_{i}^{\T} w \leq 0, \ i = 1,2,\dots,m\}$.
    \begin{enumerate}[(1)]
        \item If $S$ has bounded curvature, then $S = K + R$ for a compact set $K$.
        \item If all $f_{i}$ are convex, then $R = \operatorname{rec}(S)$, and the following are equivalent:
            \begin{enumerate}[(a)]
                \item $S$ has bounded curvature.
                \item $S = K + V$, where $K$ is a compact set and $V$ is a polyhedral cone.
                \item $S = K + \operatorname{rec}(S)$, where $K$ is a compact set.
            \end{enumerate}
    \end{enumerate}
\end{lemma}

We prove Lemma~\ref{lem:bounded_curvature} in Appendix~\ref{app:bounded}.

The following result is due to Bank and Mandel~\cite[Section 6.2, Theorem 2 with Lemma 2(ii)]{BankMandel88}.
\begin{theorem} \label{thm:bankmandel}
    Let $S = C + V$ where $C$ is compact and $V$ is a cone
    generated by finitely many integer vectors. Let $f(x) = x^{\T} D x + q^{\T} x + c$ 
    for $D \in \Q^{n \times n}$, $q \in \R^{n}$, $c \in \R$.
    If $f$ is bounded from below on $S \cap \left(\Z^{p} \times \R^{n-p}\right) \neq \emptyset$, 
    then it attains its infimum there.
\end{theorem}

This gives the following corollary, which proves one direction of Theorem~\ref{thm:main}(iii),
without convexity.
\begin{corollary} \label{cor:forward}
    Let $S$ be given by~\eqref{def:feasible_set}, where all $f_{i}$ have rational coefficients,
    and assume that $S$ has bounded curvature. Then $S$ is $\text{FW}_{n}$.
\end{corollary}

\begin{proof}
    Let $\Lambda = \tau\left(\Z^{p} \times \R^{n-p}\right)$ be a rational mixed integer lattice and
    let $P = \{x \in \R^{n} : Ax \leq b\}$ be a rational polyhedron. After the change of variables
    $x = \tau(y)$, the lattice becomes $\Z^{p} \times \R^{n-p}$, and the transformed set is still
    defined by quadratic inequalities with rational coefficients. It still has bounded curvature,
    since $\tau^{-1}(B + N) = M^{-1}(B - \bar{x}) + M^{-1} N$, and $M^{-1} N$ is the intersection of the
    kernels of the transformed quadratic forms. Hence we may assume that $\Lambda = \Z^{p} \times \R^{n-p}$.
    The set $S \cap P$ is defined by the constraints of $S$ and the linear inequalities of $P$, so it
    also has bounded curvature. By Lemma~\ref{lem:bounded_curvature}(1), we can write
    $S \cap P = C + R$ for a compact set $C$, where
    $R = \{w \in N : q_{i}^{\T} w \leq 0, \ i = 1,2,\dots,m, \ Aw \leq 0\}$.
    This cone is defined by rational linear equations and inequalities, so it is generated by finitely
    many rational vectors, which can be scaled to be integer.
    Thus Theorem~\ref{thm:bankmandel} applies to $S \cap P$ for every quadratic objective
    with rational coefficients, and $S$ is $\text{FW}_{n}$.
\end{proof}

\begin{remark} \label{rmk:nonconvex}
    The converse of Corollary~\ref{cor:forward} does not hold without convexity.
    For example, let $S := \{x \in \R^{n} : 1 - x_{1}^{2} \leq 0\}$. Its orthogonal projection onto
    $N^{\perp} = \operatorname{span}(e^{1})$ is unbounded, so $S$ does not have bounded curvature.
    However, for every rational polyhedron $P$, the set $S \cap P$ is the union of the rational polyhedra
    $P \cap \{x : x_{1} \geq 1\}$ and $P \cap \{x : x_{1} \leq -1\}$. By~\cite[Theorem 8.5]{Schrijver98},
    each of them is the sum of a polytope and a cone generated by finitely many integer vectors.
    As in the proof of Corollary~\ref{cor:forward}, Theorem~\ref{thm:bankmandel} applies to each of them.
    If $f_{0}$ is bounded from below on $S \cap P \cap \Lambda$, it attains its infimum on the intersection of
    each piece with $\Lambda$ when this is nonempty, and the smaller of these minima is the minimum over
    $S \cap P \cap \Lambda$. Hence $S$ is $\text{FW}_{n}$.
\end{remark}

We now show the converse, for sets of dimension at least 3 with a rational affine hull, 
given by convex quadratics with rational coefficients, and for $2 \leq p \leq n-1$.

The following example shows that a single convex quadratic constraint can prevent attainment.
It illustrates the construction in the proof of Theorem~\ref{thm:counterexample}.

\begin{example} \label{ex:multiple_integer_convex}
    For every $n,p \in \N$ such that $2 \leq p \leq n-1$, there 
    exists an instance of~\eqref{prob:MIQCQP} with a rational mixed integer lattice $\Lambda$ with $p$ integer
    components, a nonconvex $f_{0}$, and a single convex quadratic constraint $f_{1}$, both with rational
    coefficients, 
    such that the infimum is finite but not attained.
\end{example}

\begin{figure}[t]
    \centering
    \includegraphics[width=\textwidth]{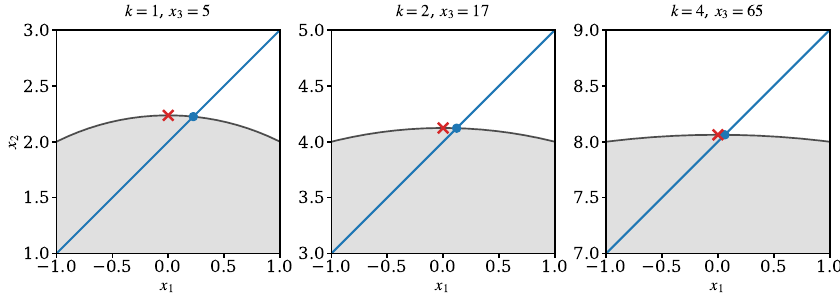}
    \caption{Example~\ref{ex:multiple_integer_convex} in the planes $x_{3} = 4k^{2} + 1$. The shaded
    region is the slice of $S$, and the line is the slice of the lattice line with $x_{2} = x_{1} + 2k$.
    The dot is the point $x^{k}$, with objective value $t_{k}^{2}$. The cross is the point with
    $x_{1} = 0$ and $x_{3} = x_{2}^{2}$, which has objective value 0 but is not in the lattice.}
    \label{fig:example2}
\end{figure}

\begin{proof}
    Consider the instance of~\eqref{prob:MIQCQP} as  
    \begin{gather*}
        \text{inf }  x_{3} - x_{2}^{2} \\ 
        \text{s.t. }  x_{1}^{2} + x_{2}^{2} - x_{3} \leq 0 \\ 
        x \in \{ (t, 2 \lambda_{1} + t, 2 \lambda_{2} + 1) : t \in \R, \lambda_{1}, \lambda_{2} \in \Z \}.
    \end{gather*}
    On the feasible region, the objective satisfies $x_{3} - x_{2}^{2} \geq x_{1}^{2} \geq 0$,
    with equality if and only if $x_{1} = 0$ and $x_{3} = x_{2}^{2}$.
    The mixed integer lattice is a union of lines that do not intersect this parabola.

    For some feasible $x = (t, 2 \lambda_{1} + t, 2 \lambda_{2}+1)$ 
    to have an objective value of 0, it would require $t=0$, but then we would 
    also need $(2 \lambda_{1})^{2} = 2 \lambda_{2} + 1$, which is impossible as 
    $4 \lambda_{1}^{2}$ is even while $2 \lambda_{2} + 1$ is odd.

    Now fix $\lambda_{1}$ and $\lambda_{2}$. The values of $t$ for which the constraint holds
    with equality are the roots of a quadratic equation, and one of them is 
    \begin{equation} \label{eq:3d_meet_constraint}
    t = \frac{-2 \lambda_{1} + \sqrt{2 - 4 \lambda_{1}^{2} + 4 \lambda_{2} }}{2}.
    \end{equation}
    If this is a real number, the corresponding point has objective value $t^{2}$.
    We show that $t$ can be made arbitrarily small. 

    For $k = 1,2,\dots$, let $\lambda_{1,k} := k$ and $\lambda_{2,k} := 2k^{2}$, so that
    $2 - 4 \lambda_{1,k}^{2} + 4 \lambda_{2,k} = 4k^{2} + 2$, and let $t_{k}$ be given
    by~\eqref{eq:3d_meet_constraint}. Then $x^{k} := (t_{k}, 2 \lambda_{1,k} + t_{k}, 2 \lambda_{2,k} + 1)$
    is feasible with objective value $t_{k}^{2}$, and
    $t_{k} = \sqrt{k^{2} + 1/2} - k = \frac{1}{2(\sqrt{k^{2} + 1/2} + k)}$.
    Hence $\lim_{k \to \infty} t_{k} = 0$, so the infimum is 0, and it is not attained.
    Figure~\ref{fig:example2} shows these points.

    By Remark~\ref{rmk:extension}, adding integer or continuous variables gives an example
    for all $n$ and $p$ with $2 \leq p \leq n-1$.
\end{proof}

\begin{theorem} \label{thm:counterexample}
    Let $S$ be given by~\eqref{def:feasible_set}, where all $f_{i}$ are convex with rational
    coefficients. Assume that $S$ has dimension at least 3, has a rational affine hull, and does not
    have bounded curvature. Then for every $p$ with $2 \leq p \leq n-1$, there exist a quadratic function $f_{0}$ with
    rational coefficients, a rational mixed integer lattice $\Lambda$ with $p$ integer components,
    and a rational polyhedron $H$ such that $f_{0}$ has a finite infimum over
    $S \cap H \cap \Lambda$, but does not attain it.
\end{theorem}

\begin{proof}[sketch]
After a rational change of variables, we may assume that $S$ is full-dimensional. Since $S$ does not
have bounded curvature, there are points $x^{k} \in S$ whose projections $u^{k}$ onto $N^{\perp}$ satisfy
$\|u^{k}\| \to \infty$. We write $x^{k} = u^{k} + w^{k}$ with $w^{k} \in N$ of norm $O(\|u^{k}\|^{2})$.
The constraints that curve in the limit direction $e$ of $u^{k}/\|u^{k}\|$ force $\|w^{k}\|$ to grow
quadratically in $\|u^{k}\|$, and the limit direction of $w^{k}/\|w^{k}\|$ is a recession direction of $S$
along which these constraints decrease. We replace these two limit directions by rational directions
$\bar{e}$ and $d$ with the same properties.

Together with a third rational direction $v$, the directions $\bar{e}$ and $d$ span a three-dimensional
affine subspace $H$ through a rational interior point of $S$, with coordinates $(s,t,y)$. For all large
$s$ and all $y$ in a small interval, points in $S \cap H$ satisfy $t \geq \psi(s,y)$, where $\psi$ is a
quadratic function with rational coefficients. The objective $f_{0} = t - \psi(s,y) + (y-\theta)^{2}$,
for a rational $\theta$, is nonnegative on $S \cap H$, and it vanishes only on the curve where $y = \theta$
and $t = \psi(s,\theta)$. The lattice in $H$ consists of the points with $s = k + \lambda \sigma$,
$t \in \Z + \frac{1}{2D}$, and $y = \sigma$, where $k \in \Z$, $\sigma \in \R$, $\lambda \in \{1,2\}$, and
$D$ is a common denominator of the values $\psi(k + \lambda\theta, \theta)$. It does not meet the curve,
since these values lie in $\frac{1}{D}\Z$. The value of $\psi$
changes at rate proportional to $k$, so a change of $\sigma$ of order $\frac{1}{k}$ reaches a point with
$t = \psi(s,y)$ and $t \in \Z + \frac{1}{2D}$. These lattice points have objective values of order
$\frac{1}{k^{2}}$, so the infimum is 0 and it is not attained. Finally, we add $p-2$ integer and $n-p-1$
continuous directions to obtain a lattice in $\R^{n}$ with $p$ integer components.
\end{proof}

Example~\ref{ex:multiple_integer_convex} follows a similar construction, with $\bar{e} = e^{2}$, $d = e^{3}$, $v = e^{1}$,
$\psi(s,y) = s^{2} + y^{2}$, and $\theta = 0$. There, the objective is $x_{3} - x_{2}^{2} = t - \psi(s,y) + y^{2}$,
and the lattice consists of the points with $s - y \in 2\Z$ and $t \in 2\Z + 1$.
The full proof of Theorem~\ref{thm:counterexample} is in Appendix~\ref{app:bounded}.

Together, these results give part~(iii) of Theorem~\ref{thm:main}.

\begin{corollary} \label{cor:characterization}
    Let $S$ be given by~\eqref{def:feasible_set}, where all $f_{i}$ are convex with rational coefficients.
    Assume that $S$ has dimension at least 3 and a rational affine hull, and let $2 \leq p \leq n$.
    Then $S$ is $\text{FW}_{p}$ if and only if $S$ has bounded curvature.
\end{corollary}

\begin{proof}
    By Corollary~\ref{cor:forward}, if $S$ has bounded curvature, then it is $\text{FW}_{n}$,
    and hence $\text{FW}_{p}$. Conversely, if $S$ does not have bounded
    curvature, Theorem~\ref{thm:counterexample} gives a counterexample with $2$ integer components,
    which exists since $n \geq 3$. Hence $S$ is not $\text{FW}_{p}$ for any $p \geq 2$.
\end{proof}

Every full-dimensional set with $n \geq 3$ satisfies the assumptions of Corollary~\ref{cor:characterization}.
So does every set of dimension at least 3 defined by a single convex quadratic inequality and linear
inequalities with rational coefficients, since its affine hull is rational by~\cite[Theorem 1.2]{DelPia24}.
For a single convex quadratic inequality, we obtain the following characterization, which we prove in
Appendix~\ref{app:bounded}.
\begin{corollary} \label{cor:one_constraint}
    Let $n \geq 3$, and let $S := \{x \in \R^{n} : x^{\T} Q x + q^{\T} x + c \leq 0\}$, where $Q$ is
    positive semidefinite and nonzero, and $Q$, $q$, and $c$ are rational. Then $S$ is $\text{FW}_{2}$
    if and only if $q \in \operatorname{range} Q$. In this case, $S$ is $\text{FW}_{n}$.
\end{corollary}

To complete the picture, we also have the following elementary result~\cite[Section 7.3, Remark 4]{BankMandel88}.
\begin{theorem} \label{thm:all_integer}
    Let $X \subseteq \Z^{n}$ be nonempty. If $f$ is a polynomial with rational coefficients 
    which is bounded from below on $X$, then $f$ attains its infimum on $X$.
\end{theorem}
This gives the following corollary for mixed integer lattices with $n$ integer components.
\begin{corollary} \label{cor:all_integer}
    Let $S$ be given by~\eqref{def:feasible_set} with $n \geq 1$, let $\Lambda$ be a rational mixed integer
    lattice with $n$ integer components, let $P$ be a rational polyhedron, and let $f_{0}$ be a quadratic
    function with rational coefficients that is bounded from below on $S \cap P \cap \Lambda \neq \emptyset$.
    Then $f_{0}$ attains its infimum on $S \cap P \cap \Lambda$. Consequently, $S$ is $\text{FW}_{n}$ if and
    only if $S$ is $\text{FW}_{n-1}$, and if $n \leq 2$, then $S$ is $\text{FW}_{n}$ if and only if $S$ is
    $\text{FW}_{0}$.
\end{corollary}

\begin{proof}
    Write $\Lambda = \tau(\Z^{n})$. Then $f_{0} \circ \tau$ is a polynomial with rational coefficients that is
    bounded from below on the nonempty set $\tau^{-1}(S \cap P \cap \Lambda) \subseteq \Z^{n}$.
    By Theorem~\ref{thm:all_integer}, it attains its infimum there, so $f_{0}$ attains its infimum on
    $S \cap P \cap \Lambda$. By definition, every $\text{FW}_{n}$ set is $\text{FW}_{n-1}$, and the converse
    follows from the first statement. For $n \leq 2$, the last statement then follows from
    Corollary~\ref{cor:quadratic_one_integer} in Section~\ref{sec:one_integer}.
\end{proof}

Finally, we consider sets of dimension at most 2. This proves part~(ii) of Theorem~\ref{thm:main},
without convexity.
\begin{corollary} \label{cor:low_dimension}
    Let $S$ be given by~\eqref{def:feasible_set}, and assume that the affine hull of $S$ is rational
    and has dimension at most 2. Then $S$ is $\text{FW}_{n}$ if and only if $S$ is $\text{FW}_{0}$.
\end{corollary}

\begin{proof}[sketch]
By~\cite[Theorem 1.1]{DelPia24}, a rational affine map identifies the affine hull of $S$, intersected with
$\Lambda$, with a mixed integer lattice in $\R^{k}$, where $k \leq 2$. The corresponding set in $\R^{k}$ is
$\text{FW}_{0}$, and hence $\text{FW}_{k}$ by Corollary~\ref{cor:all_integer}. The full proof is in
Appendix~\ref{app:bounded}.
\end{proof}

\section{Conditions on the objective}
\label{sec:convex}

In this section, we show that the infimum of~\eqref{prob:MIQCQP} is attained under
assumptions on the objective that are not covered by Theorem~\ref{thm:main}.
The main result is Theorem~\ref{thm:convex}: if the data are rational, the constraints are
convex, and $f_{0}$ is bounded from above and from below on $S \cap \Lambda$, then the
infimum is attained.

The idea of the proof of the following theorem is as follows.
If $S \cap \Lambda$ is bounded, then the infimum is attained. Otherwise, there is a rational recession
direction $d$ of the constraints. Since $f_{0}$ is bounded on $S \cap \Lambda$, it is constant along $d$
on $S \cap \Lambda$. After adding two linear constraints, it is constant along $d$ on the continuous
relaxation, and we can project out $d$ and use induction on $n$.

\begin{theorem} \label{thm:convex}
    Consider an instance of~\eqref{prob:MIQCQP} with rational data.
    If all constraints $f_{i}$ $(i=1,2,\dots,m)$ are convex quadratic functions, and $f_{0}(x)$ is bounded from above on $S \cap \Lambda$, 
    then~\eqref{prob:MIQCQP} attains its infimum when it is finite. 
\end{theorem}

We prove Theorem~\ref{thm:convex} in Appendix~\ref{app:convex}, using three lemmas.

\begin{corollary} \label{cor:quasi-convex}
    Consider an instance of~\eqref{prob:MIQCQP} with rational data, where
    all constraints $f_{i}$ $(i=1,\dots,m)$ are convex. Assume that $f_{0}$ is convex,
    or that $f_{0}$ is quasi-convex over a full-dimensional convex set $C \supseteq S \cap \Lambda$.
    Then~\eqref{prob:MIQCQP} attains its infimum (when it is finite).
\end{corollary}

\begin{proof}
    First, assume that $f_{0}$ is convex, that is, $Q_{0}$ is positive semidefinite.
    Let $r \in \Q$ with $r > \inf_{y \in S \cap \Lambda} f_{0}(y)$,
    and add the constraint $f_{0}(x) - r \leq 0$.
    This constraint is convex and has rational coefficients.
    A point attains the infimum of the original problem if and only if it attains
    the infimum of the new problem. In the new problem, $f_{0}$ is bounded from above
    on the feasible region, so the infimum is attained by Theorem~\ref{thm:convex}.

    Next, assume that $Q_{0}$ is not positive semidefinite. Then $f_{0}$ is quasi-convex but not
    convex over a full-dimensional convex set $C \supseteq S \cap \Lambda$. By~\cite[Theorem 6.3]{GeneralizedConcavity10},
    applied to $-f_{0}$, the function $f_{0}$ is bounded from above on $C$ (see also~\cite[Lemma 7]{LuoZhang}),
    and thus on $S \cap \Lambda$. The infimum is attained by Theorem~\ref{thm:convex}.
\end{proof}

If $S$ is full-dimensional and $f_{0}$ is quasi-convex over $S$, we can take $C = S$.
If the polyhedron defined by the linear constraints is full-dimensional and $f_{0}$ is quasi-convex over it,
as in~\cite[Theorem 3]{LuoZhang}, we can take $C$ to be this polyhedron.

For at most one integer variable, Corollary~\ref{cor:one_integer_quasi-convex} in
Appendix~\ref{app:convex} removes the rationality assumption of Corollary~\ref{cor:quasi-convex}, and it
only requires that the objective is quasi-convex over $S$.

\begin{remark} \label{rmk:tightness}
    Theorem~\ref{thm:convex} and Corollary~\ref{cor:quasi-convex} restrict the objectives in
    Theorem~\ref{thm:counterexample}. There, the constraints that define $S \cap H$ are convex with
    rational coefficients. Hence the objective $f_{0}$ is not convex, it is not quasi-convex over any
    full-dimensional convex set containing $S \cap H$, and it is unbounded from above on $S \cap H \cap \Lambda$.
    In particular, in this setting every linear objective with rational coefficients attains its infimum
    when it is finite. The objective $x_{3} - x_{2}^{2}$ of Example~\ref{ex:multiple_integer_convex}
    has these properties.
\end{remark}

\paragraph{AI Declaration.}
AI tools (Claude, developed by Anthropic) were used for proofreading and
suggesting improvements to the exposition. We take full responsibility for the content of the paper.

\paragraph{Funding.}
The authors are partially supported by AFOSR Grant \mbox{FA9550-23-1-0433}
and ONR Grant \mbox{N00014-25-1-2490}. Any opinions, findings, conclusions,
or recommendations expressed in this material are those of the authors and
do not necessarily reflect the views of the Air Force Office of Scientific
Research or the Office of Naval Research.

%
% ---- Bibliography ----
%
\bibliographystyle{plain}
\bibliography{references}

\clearpage
\appendix

\section{Auxiliary lemmas}
\label{app:prelim}

We collect some lemmas that are used in several proofs.
The first lemma describes the recession directions of a set defined by convex quadratic inequalities.
\begin{lemma} \label{lem:convex_recession}
    Let $X := \{x \in \R^{n} : g_{i}(x) \leq 0, \ i = 1,2,\dots,m\}$, where each
    $g_{i}(x) = x^{\T} Q_{i} x + q_{i}^{\T} x + c_{i}$ is convex, and assume that $X$ is nonempty.
    Then $\operatorname{rec}(X) = \{d \in \R^{n} : Q_{i} d = 0, \ q_{i}^{\T} d \leq 0, \ i = 1,2,\dots,m\}$.
    Moreover, $d \in \operatorname{rec}(X)$ if and only if $d$ is a recession direction of each
    inequality $g_{i}(x) \leq 0$.
\end{lemma}

\begin{proof}
    Let $d$ satisfy $Q_{i} d = 0$ and $q_{i}^{\T} d \leq 0$ for every $i$. Then
    $g_{i}(x + td) = g_{i}(x) + t q_{i}^{\T} d \leq g_{i}(x)$ for every $x \in \R^{n}$, $t \geq 0$, and $i$.
    Hence $d$ is a recession direction of each inequality, and of $X$.
    Conversely, let $x \in X$, and let $d$ be a recession direction of $X$ or of each inequality.
    In both cases, for every $i$ and every $t \geq 0$ we have
    \[
        0 \geq g_{i}(x + td) = g_{i}(x) + t\left(2 x^{\T} Q_{i} d + q_{i}^{\T} d\right) + t^{2} d^{\T} Q_{i} d .
    \]
    Dividing by $t^{2}$ and letting $t \to \infty$ gives $d^{\T} Q_{i} d \leq 0$. Since $Q_{i}$ is positive
    semidefinite, this implies $Q_{i} d = 0$. Then $g_{i}(x) + t q_{i}^{\T} d \leq 0$ for every $t \geq 0$,
    so $q_{i}^{\T} d \leq 0$.
\end{proof}

The next lemma replaces a nonzero vector in a rational polyhedron by a nearby nonzero rational vector.

\begin{lemma} \label{lem:rational_solution}
    Let $P := \{x \in \R^{n} : Ax \leq b\}$ with $A \in \Q^{m \times n}$ and $b \in \Q^{m}$,
    and let $\bar{x} \in P \setminus \{0\}$. Then for every $\epsilon > 0$ there exists
    $x \in \left(P \setminus \{0\}\right) \cap \Q^{n}$ with $\|x - \bar{x}\| < \epsilon$.
\end{lemma}

\begin{proof}
    Let $B$ be a box defined by rational inequalities, with side length less than
    $\frac{\epsilon}{\sqrt{n}}$, such that $\bar{x} \in B$ and $0 \notin B$.
    Then $P \cap B$ is a nonempty polyhedron defined by rational inequalities, since it contains $\bar{x}$.
    By~\cite[Theorem 10.1]{Schrijver98}, it contains a rational point $x$.
    This point is nonzero since $0 \notin B$, and $\|x - \bar{x}\| < \epsilon$ since both points lie in $B$.
\end{proof}

The following lemma is implicit in the proof of~\cite[Theorem 1]{LuoZhang}.
\begin{lemma} \label{lem:limit}
    Consider a convex quadratic inequality
    \begin{equation}
    \label{eq:convex_inequality}
    f(x) = x^{\T} Q x + q^{\T} x + c \leq 0.
    \end{equation}
    Let $\left( x^{k} \right)$ be a sequence of solutions of~\eqref{eq:convex_inequality}
    with $\lim_{k \to \infty} \|x^{k}\| = \infty$ such that the limit
    $d := \lim_{k \to \infty} \frac{x^{k}}{\|x^{k}\|}$ exists.
    Then $d$ is a recession direction of~\eqref{eq:convex_inequality}.  
\end{lemma}

\begin{proof}
    Taking the inequality for $x^{k}$ and dividing through by $\|x^{k}\|^{2}$ yields 
    \[\left(\frac{x^{k}}{\|x^{k}\|}\right)^{\T} Q \left(\frac{x^{k}}{\|x^{k}\|}\right) 
    + q^{\T} \frac{x^{k}}{\|x^{k}\|^{2}} + \frac{c}{\|x^{k}\|^{2}} \leq 0 . \]
    Taking the lim sup as $k$ tends to infinity gives $d^{\T} Q d \leq 0$, but since $Q$ is positive 
    semidefinite this means $Qd = 0$. 

    If we instead only divide by $\|x^{k}\|$ we get 
    \[\frac{1}{\|x^{k}\|} (x^{k})^{\T} Q x^{k} + q^{\T} \frac{x^{k}}{\|x^{k}\|} + \frac{c}{\|x^{k}\|} \leq 0 .\]
    Then we must have $q^{\T} d \leq 0$, because the first term is always nonnegative and the last term 
    tends to 0, so $q^{\T} d > 0$ would imply $f(x^{k}) > 0$ for large enough values of $k$. 
    
    Then Lemma~\ref{lem:convex_recession} gives that $d$ is a recession direction for~\eqref{eq:convex_inequality}.
\end{proof}

The following lemma is used several times. It shows that moving along a rational direction
by integer multiples of a suitable step stays in a mixed integer lattice.

\begin{lemma} \label{lem:mixed_integer_spacing}
    Let $d \in \Q^{n}$, let $p \in \{0, 1,\dots,n\}$, and let $\Lambda = \tau\left(\Z^{p} \times \R^{n-p}\right)$
    for an invertible rational affine map $\tau$. Then there exists a rational $t > 0$ such that
    $x + ktd \in \Lambda$ for every $x \in \Lambda$ and every $k \in \Z$.
\end{lemma}

\begin{proof}
    After the change of variables $x = \tau(y)$, we may assume that $\Lambda = \Z^{p} \times \R^{n-p}$,
    since $d$ is replaced by the rational vector $M^{-1} d$.
    Write $d_{i} = \frac{a_{i}}{b_{i}}$ with $a_{i} \in \Z$ and $b_{i} \in \Z_{> 0}$.
    Then $t := \prod_{i=1}^{p} b_{i}$, where the empty product is 1, ensures that the
    first $p$ components of $x + ktd$ are integer.
\end{proof}

\section{One integer variable}
\label{app:one_integer}

This appendix contains the proof of Theorem~\ref{thm:generic_one_integer} and further consequences of it.

\begin{proof}[of Theorem~\ref{thm:generic_one_integer}]
    After the change of variables $x = My + \bar{x}$, we may assume that $\Lambda = \Z \times \R^{n-1}$
    and that the relaxations in the assumption are those with bounds $L,U \in \Z \cup \{-\infty, \infty\}$
    on $x_{1}$. Indeed, $(a^{*})^{\T} x = y_{1} + (a^{*})^{\T}\bar{x}$.
    To refer to the feasible regions of the problems, we define
    \[K := \{x \in \R^{n} : g_{i}(x) \leq 0, \quad i =1,2,\dots, m \},\]
    \[K_{L,U} := \{x \in \R^{n} : L \leq x_{1} \leq U, g_{i}(x) \leq 0, \quad i=1,2,\dots,m \},\]
    and $T := \{x \in \R^{n} : g_{0}(x) \leq c \}$. By assumption, $T \cap K$ has finitely many connected components.

    Now, assume for contradiction that~\eqref{prob:generic_one_integer} does not attain its infimum.
    Then no $x \in K \cap (\Z \times \R^{n-1})$ satisfies $g_{0}(x) \leq c$, and there is a sequence $(x^{k})$
    in $K \cap (\Z \times \R^{n-1})$ with $\lim_{k \to \infty} g_{0}(x^{k}) = c$.

    First, assume that the sequence $(x^{k}_{1})$ is bounded. Since its elements are integers, it takes
    finitely many values $j$. For each such $j$, the relaxation~\eqref{prob:relaxed_with_linear} with $L = U = j$
    has feasible region $K_{j,j} \subseteq K \cap (\Z \times \R^{n-1})$, so its infimum $\gamma_{j}$ is finite and
    at least $c$. By assumption, $\gamma_{j}$ is attained, so $\gamma_{j} > c$. Then $g_{0}(x^{k}) \geq \min_{j} \gamma_{j} > c$
    for every $k$, a contradiction.

    Hence $(x^{k}_{1})$ is unbounded. Passing to a subsequence, we may assume that $x^{k}_{1}$
    tends to $\infty$, as the case where it tends to $-\infty$ is symmetric.
    Every $x \in T \cap K$ satisfies $x_{1} \notin \Z$, since otherwise $x$ would attain the infimum.
    Let $C$ be a connected component of $T \cap K$. Its projection onto the first coordinate is connected,
    and thus an interval in $\R$ that contains no integer. Hence there is $n_{C} \in \Z$ such that
    $n_{C} < x_{1} < n_{C}+1$ for every $x \in C$.
    Since $T \cap K$ has finitely many connected components, we can take $B := \max_{C} n_{C} + 1$,
    or any $B \in \Z$ if $T \cap K$ is empty. Then $T \cap K_{B, \infty}$ is empty. Hence every feasible point
    of the relaxation~\eqref{prob:relaxed_with_linear} with $L = B$ and $U = \infty$ has objective value
    larger than $c$.
    On the other hand, $x^{k} \in K_{B, \infty}$ for all $k$ large enough, and $\lim_{k \to \infty} g_{0}(x^{k}) = c$.
    Hence the infimum of this relaxation is $c$. By assumption, it is attained, so there is $x \in K_{B,\infty}$
    with $g_{0}(x) = c$. This contradicts that $T \cap K_{B, \infty}$ is empty.
\end{proof}

Theorem~\ref{thm:generic_one_integer} also lifts other known attainment results from
continuous variables to one integer variable. This applies to every class of problems defined by
polynomials that remains in the class when linear inequalities are added. The sets in
Theorem~\ref{thm:generic_one_integer} are then semialgebraic, so they have finitely many connected
components. If the continuous result allows linear inequalities with real coefficients, the mixed integer
lattice $\Lambda$ need not be rational. The proofs are the same as the proof of
Corollary~\ref{cor:quadratic_one_integer}, so we omit them.

For example, we can generalize~\cite[Theorem 2]{LuoZhang}.
\begin{corollary} \label{cor:one_integer_convex}
    Consider a problem of the form~\eqref{prob:generic_one_integer}, where $g_{0}$ is a quadratic function,
    all constraint functions $g_{i}$ $(i=1,2,\dots,m)$ are convex quadratic functions,
    and at most one of them is not affine.
    Then~\eqref{prob:generic_one_integer} attains its infimum (when it is finite).
\end{corollary}

The following lifts the cubic programming attainability from Klatte~\cite{Klatte19} 
to the mixed integer setting.
\begin{corollary} \label{cor:cubic_one_integer}
    Consider a problem of the form 
    \begin{align}
    \label{prob:cubic}
    \begin{split}
        \inf\; &g_{0}(x) \\
        \text{s.t.}\; & Ax \leq b \\ 
        & x \in \Lambda,
    \end{split}
    \end{align}
    where $A \in \R^{m \times n}$, $b \in \R^{m}$, $g_{0}$ is a polynomial of degree at most 3,
    and $\Lambda$ is as in~\eqref{prob:generic_one_integer}.

    Then~\eqref{prob:cubic} attains its infimum (when it is finite).
\end{corollary}

We can also lift the result of Belousov and Klatte for convex polynomials~\cite[Theorem 3]{BelousovKlatte02}.
\begin{corollary} \label{cor:generic_convex}
    Consider a problem of the form 
    \begin{align}
    \label{prob:convex}
    \begin{split}
        \inf\; &g_{0}(x) \\ 
        \text{s.t.}\; & g_{i}(x) \leq 0, \quad i=1,2,\dots,m \\
        & x \in \Lambda,
    \end{split}
    \end{align}
    where $g_{j} : \R^{n} \to \R$ for $j=0,1,\dots,m$ is a convex polynomial,
    and $\Lambda$ is as in~\eqref{prob:generic_one_integer}.

    Then~\eqref{prob:convex} attains its infimum (when it is finite).
\end{corollary}

Corollaries~\ref{cor:one_integer_convex}, \ref{cor:cubic_one_integer}, and~\ref{cor:generic_convex}
do not require rational data or a rational mixed integer lattice. In contrast, the next section shows
that rational data are needed with two integer variables.

\section{Irrational data}
\label{sec:irrational}

In this section, we show that rational data are needed in general once there are two integer
variables, even for linear constraints.

We first recall two results of Meyer. Theorem~\ref{thm:irrational_linear} is the example
in~\cite[Section 2]{Meyer74}, and Theorem~\ref{thm:rational_constraints} follows
from~\cite[Theorem 2.1]{Meyer74}.

\begin{theorem} \label{thm:irrational_linear}
    There exists an instance of~\eqref{prob:MIQCQP} with $n=p=2$ and $\Lambda$ rational,
    where all constraints and the objective function are affine, 
    but some constraint has irrational coefficients, such that the infimum is 
    finite but not attained.
\end{theorem}

\begin{theorem} \label{thm:rational_constraints}
    Consider an instance of~\eqref{prob:MIQCQP} where $\Lambda$ is rational, all functions
    $f_{j}$ $(j=0,\dots,m)$ are affine, and all constraints $f_{i}$ $(i=1,\dots,m)$
    have rational coefficients. Then the infimum is attained whenever it is finite.
\end{theorem}

Hence, when all functions are affine, rationality of the constraints alone ensures
that the infimum is attained, and this assumption cannot be dropped in general.
The next example shows that rational constraints are not enough for a quadratic objective:
an objective with irrational coefficients may not attain its infimum.

\begin{example} \label{ex:irrational_objective}
    There exist instances of~\eqref{prob:MIQCQP} where $\Lambda$ is rational, all constraints
    $f_{i}$ $(i=1,\dots,m)$ are affine with rational coefficients, 
    but the objective $f_{0}$ is a quadratic function with irrational coefficients, 
    such that the infimum is finite but not attained.
\end{example}

\begin{proof} 
    Consider
    \begin{align*} 
        \text{inf } & (x_{1} - x_{2} \sqrt{2})^{2} \\
        \text{s.t. } & x_{2} \geq 1 \\ 
        & x \in \Z^{2}. 
    \end{align*}
    An objective value of 0 cannot be attained for any $x \in \Z^{2}$, because that would require $\frac{x_{1}}{x_{2}} = \sqrt{2}$.
    For $k=1,2,\dots$, define positive integers $a_{k}, b_{k}$ by $a_{k} + b_{k} \sqrt{2} := (1 + \sqrt{2})^{k}$,
    and let $x^{k} := (a_{k}, b_{k})$. Then $a_{k} - b_{k} \sqrt{2} = (1 - \sqrt{2})^{k}$, so the objective
    value of $x^{k}$ is $(\sqrt{2} - 1)^{2k}$, which tends to 0.
    Since $b_{k} \geq 1$, the points $x^{k}$ are feasible. Hence the infimum is 0, and it is not attained.
\end{proof}

Both examples have two integer variables. With one integer variable, the infimum is attained
in both settings, even if the mixed integer lattice is not rational. This follows from
Corollary~\ref{cor:cubic_one_integer}, since the objectives are polynomials of degree at most 2
and the constraints are linear.

Table~\ref{tab:summary} summarizes when integer variables can destroy attainment for a set that is
$\text{FW}_{0}$, with and without rational data.

\begin{table}[t]
\caption{When can integer variables destroy attainment for a set $S$ given by~\eqref{def:feasible_set} that is
$\text{FW}_{0}$, convex or not? Here $p$ is the number of integer components of the mixed integer lattice.}
\label{tab:summary}
\centering
\small
\begin{tabular}{p{2.2cm}p{4.6cm}p{4.3cm}}
\hline
 & Rational data & Irrational data \\
\hline
$p \leq 1$ & never (Corollary~\ref{cor:quadratic_one_integer}) & never, if the continuous result allows real data (Section~\ref{sec:one_integer}) \\
$2 \leq p \leq n-1$ & only if $S$ does not have bounded curvature (Corollary~\ref{cor:forward}), and exactly then if $S$ is convex with a rational affine hull of dimension at least 3 (Theorem~\ref{thm:main}(iii)) & possible, even for polyhedra (Theorem~\ref{thm:irrational_linear} and Remark~\ref{rmk:extension}) \\
$p = n \geq 2$ & never (Corollary~\ref{cor:all_integer}) & possible, even for polyhedra (Theorem~\ref{thm:irrational_linear} and Remark~\ref{rmk:extension}) \\
\hline
\end{tabular}
\end{table}

\section{Two or more integer variables}
\label{app:bounded}

This appendix contains the proofs omitted from Section~\ref{sec:bounded_curvature}.

\begin{proof}[of Lemma~\ref{lem:bounded_curvature}]
    Let $U := N^{\perp}$, and let $\pi_{U}$ and $\pi_{N}$ denote the orthogonal projections onto $U$ and $N$.
    For every $x \in \R^{n}$ and $w \in N$, we have $f_{i}(x + w) = f_{i}(x) + q_{i}^{\T} w$, since
    $Q_{i} w = 0$. In particular, $x + r \in S$ for every $x \in S$ and $r \in R$.

    For (1), write each $x \in S$ as $x = u + w$, where $u = \pi_{U}(x)$ and $w = \pi_{N}(x)$.
    Since $S$ has bounded curvature, $\pi_{U}(S)$ is contained in a bounded set $B$.
    For each $u \in \pi_{U}(S)$, the set of $w \in N$ such that $u+w \in S$ is
    \[F_{u} := \{w \in N : q_{i}^{\T} w \leq -f_{i}(u), \quad i=1,2,\dots,m \}. \]
    This polyhedron is nonempty, and its recession cone is $R$.
    By~\cite[Theorem 8.5]{Schrijver98}, we can write $F_{u} = \operatorname{conv}(Z_{u}) + R$,
    where $Z_{u} \subseteq F_{u}$ is a finite set that contains an arbitrary point of each minimal face of $F_{u}$.

    Each minimal face of $F_{u}$ is of the form $\{w \in N : q_{i}^{\T} w = -f_{i}(u), \ i \in I'\}$
    for some $I' \subseteq \{1,2,\dots,m\}$, and we choose its point of minimum norm. This point is
    $G_{I'} b_{I'}(u)$, where $b_{I'}(u) := (-f_{i}(u))_{i \in I'}$ and $G_{I'}$ is a matrix that depends
    only on $I'$ and on the constraints. Since $u \in B$ is bounded, $f_{i}(u)$ is also bounded.
    Hence each element of $Z_{u}$ lies in a bounded set $B'$ that is independent of $u$.
    We can take $B'$ convex, so that $\operatorname{conv}(Z_{u}) \subseteq B'$. Then we can write
    $w = \bar{w} + r$, where $\bar{w} \in B' \cap F_{u}$ and $r \in R$. Since $u$ and $\bar{w}$ are bounded,
    we have $\|u + \bar{w}\| \leq \rho$ for a constant $\rho$, and $u + \bar{w} \in S$ since $\bar{w} \in F_{u}$.
    Hence $S \subseteq K + R$, where $K := S \cap \{x \in \R^{n} : \|x\| \leq \rho\}$ is compact.
    The reverse inclusion $K + R \subseteq S$ holds by the first paragraph of the proof.

    For (2), Lemma~\ref{lem:convex_recession} gives $\operatorname{rec}(S) = R$. Then (a) $\Rightarrow$ (c)
    follows from (1), and (c) $\Rightarrow$ (b) holds since $R$ is a polyhedral cone.
    For (b) $\Rightarrow$ (a), we have $V \subseteq \operatorname{rec}(S) = R \subseteq N$.
    Hence $\pi_{U}(S) \subseteq \pi_{U}(K)$, which is bounded.
\end{proof}

The next lemma characterizes bounded curvature for a set defined by a single convex quadratic inequality.

\begin{lemma} \label{lem:one_constraint}
    Let $S$ be given as in~\eqref{def:feasible_set} for a single 
    convex quadratic constraint $f_{1}$ with $Q_{1} \neq 0$.
    Then $S$ has bounded curvature if and only if $q_{1} \in \operatorname{range} Q_{1}$. 
    Moreover, if $q_{1} \notin \operatorname{range} Q_{1}$ then $S$ is full-dimensional.
\end{lemma}

\begin{proof}
    Let $N = \ker Q_{1}$ and $U = N^{\perp}$, 
    and let $\pi_{U}$ denote orthogonal projection onto $U$.
    If $q_{1} \in \operatorname{range} Q_{1}$, 
    let $\bar{x} \in \R^{n}$ satisfy $q_{1} = -2Q_{1}\bar{x}$. Then for any $x \in \R^{n}$ 
    we can write
    \[x^{\T} Q_{1} x + q_{1}^{\T} x + c_{1} = (x-\bar{x})^{\T} Q_{1} (x-\bar{x}) + c' , \]
    for some constant $c'$. Then $(x-\bar{x})^{\T} Q_{1} (x-\bar{x}) = u^{\T} Q_{1} u$, where 
    $u = \pi_{U}(x-\bar{x})$. But $Q_{1}$ is positive 
    definite on $U$, so the set defined by $u^{\T} Q_{1} u + c' \leq 0$ is an ellipsoid, 
    which is bounded. Thus, $\pi_{U}(S)$ is bounded, so $S$ has bounded curvature.

    Now suppose $q_{1} \notin \operatorname{range} Q_{1}$, so we can 
    write $q_{1} = q_{U} + q_{N}$ for $q_{U} \in U, q_{N} \in N \setminus \{0\}$. 
    Then we have $Q_{1} q_{N} = 0$ and $q_{1}^{\T} q_{N} = \|q_{N}\|^{2}$, 
    so for $u \in U, t \in \R$ we can write 
    \[f_{1}(u - t q_{N}) = u^{\T} Q_{1} u + q_{U}^{\T} u + c_{1} - t \|q_{N}\|^{2} .\]
    For large enough $t$, this is negative, so every $u \in U$ belongs to $\pi_{U}(S)$. 
    Since $Q_{1} \neq 0$, $U$ is unbounded, showing $S$ does not have bounded curvature. 
    Moreover, this also shows there are $x \in S$ with $f_{1}(x) < 0$, which is an interior 
    point since $f_{1}$ is continuous.
\end{proof}

\begin{proof}[of Corollary~\ref{cor:one_constraint}]
    If $q \in \operatorname{range} Q$, then $S$ has bounded curvature by Lemma~\ref{lem:one_constraint},
    so $S$ is $\text{FW}_{n}$ by Corollary~\ref{cor:forward}. Otherwise, $S$ is full-dimensional and
    does not have bounded curvature by Lemma~\ref{lem:one_constraint}, so $S$ is not $\text{FW}_{2}$ by
    Corollary~\ref{cor:characterization}.
\end{proof}

\begin{proof}[of Theorem~\ref{thm:counterexample}]
    We first reduce to the case where $S$ is full-dimensional. Let $\ell \geq 3$ be the dimension of $S$.
    Since the affine hull of $S$ is rational, it is $\{x_{0} + Wz : z \in \R^{\ell}\}$ for some $x_{0} \in \Q^{n}$
    and $W \in \Q^{n \times \ell}$ with full column rank. The set
    $S_{W} := \{z \in \R^{\ell} : f_{i}(x_{0} + Wz) \leq 0, \ i = 1,2,\dots,m\}$
    is full-dimensional and defined by convex quadratic inequalities with rational coefficients,
    and $S = x_{0} + W S_{W}$. If $S_{W} = K_{W} + \operatorname{rec}(S_{W})$ for a compact set $K_{W}$,
    then $S = (x_{0} + W K_{W}) + \operatorname{rec}(S)$, since $\operatorname{rec}(S) = W \operatorname{rec}(S_{W})$.
    Hence, by Lemma~\ref{lem:bounded_curvature}, $S_{W}$ does not have bounded curvature.
    Let $p_{W}$ be an integer with $\max\{2, p - n + \ell\} \leq p_{W} \leq \min\{\ell-1, p\}$,
    which exists since $\ell \geq 3$ and $2 \leq p \leq n-1$.
    Assume that the theorem holds for $S_{W}$ and $p_{W}$, with objective $f_{0}^{W}$,
    rational mixed integer lattice $\Lambda_{W} \subseteq \R^{\ell}$, and rational polyhedron $H_{W}$.
    Let $g^{1}, \dots, g^{n-\ell} \in \Q^{n}$ be such that the columns of $W$ and these vectors form a basis
    of $\R^{n}$, and define
    \[
        \Lambda := \Big\{x_{0} + Wz + \sum_{j=1}^{n-\ell} \lambda_{j} g^{j}
        : z \in \Lambda_{W}, \ \lambda \in \Z^{p - p_{W}} \times \R^{n - \ell - p + p_{W}} \Big\}.
    \]
    This is a rational mixed integer lattice with $p$ integer components, and
    $\Lambda \cap \{x_{0} + Wz : z \in \R^{\ell}\} = x_{0} + W\Lambda_{W}$, since the vectors $g^{j}$ and the columns
    of $W$ are linearly independent. Let $W^{+} := (W^{\T} W)^{-1} W^{\T}$, let $H := x_{0} + W H_{W}$,
    which is a rational polyhedron, and let $f_{0}(x) := f_{0}^{W}(W^{+}(x - x_{0}))$, which is a quadratic
    function with rational coefficients. Since $W^{+} W$ is the identity, the map $z \mapsto x_{0} + Wz$ is a
    bijection from $S_{W} \cap H_{W} \cap \Lambda_{W}$ to $S \cap H \cap \Lambda$ that preserves objective values.
    Hence $f_{0}$ has a finite infimum over $S \cap H \cap \Lambda$ that is not attained.
    Therefore, from now on we assume that $S$ is full-dimensional, so that $n \geq 3$.

    Let $N = \bigcap_{i=1}^{m} \ker Q_{i}$ and $U = N^{\perp}$.
    Our goal is to find directions in the set $S$ with specific properties, which we can use 
    to build our counterexample. 

    Let $\pi_{U}$ denote the orthogonal projection map onto $U$, and $\pi_{N}$ likewise for $N$.
    Then $S$ not having bounded curvature implies there is a sequence $(x^{k})$ such that 
    $\lim_{k \to \infty} \|\pi_{U}(x^{k})\| = \infty$. Let $u^{k} = \pi_{U}(x^{k})$. 
    Then if we consider 
    \[F_{k} := \{w \in N: q_{i}^{\T} w \leq -f_{i}(u^{k}) \quad i=1,2,\dots,m \}, \]
    we have $u^{k} + w \in S$ iff $w \in F_{k}$, and $F_{k}$ is nonempty since it contains $\pi_{N}(x^{k})$.
    First, we note that for each constraint, there is a constant $L_{i}$ such that for every $u \in \R^{n}$, 
    $|f_{i}(u)| \leq L_{i}(1 + \|u\|^{2})$. Combined with Hoffman's error bound~\cite{Hoffman52} applied 
    to 0, we get that there is a constant $L$ independent of $k$ such that there is $w^{k} \in F_{k}$ 
    with $\|w^{k}\| \leq L(1 + \|u^{k}\|^{2})$. Then we replace $x^{k}$ with $u^{k} + w^{k}$, 
    and work with this sequence instead.

    Now, pass to a subsequence such that the limit
    $\lim_{k \to \infty} \left( \frac{u^{k}}{\|u^{k}\|} \right) = e \in U$ exists.
    We show below that $w^{k} \neq 0$ for all sufficiently large $k$. Passing to a further subsequence,
    we may assume that the limit $\lim_{k \to \infty} \left( \frac{w^{k}}{\|w^{k}\|} \right) = \tilde{d}$ also exists.
    We show that $e$ and $\tilde{d}$ have some useful properties. First, since $e \in U \setminus \{0\}$, 
    we have $e \notin N$, so the set 
    \[I := \{i \in \{1,2,\dots,m\} : e^{\T} Q_{i} e > 0 \},\]
    is nonempty. For $\tilde{d}$, let $i \in I$ and let $a_{i} = e^{\T} Q_{i} e > 0$. 
    Then we have that $\lim_{k \to \infty} \frac{f_{i}(u^{k})}{\|u^{k}\|^{2}} = a_{i}$, 
    so for large enough $k$ we have $f_{i}(u^{k}) \geq \frac{1}{2} a_{i} \|u^{k}\|^{2}$. 
    Then we have $q_{i}^{\T} w^{k} \leq -f_{i}(u^{k}) \leq -\frac{1}{2} a_{i} \|u^{k}\|^{2}$, 
    so $q_{i} \neq 0$ and $\|w^{k}\| \geq \frac{a_{i}}{2 \|q_{i}\|} \|u^{k}\|^{2}$, 
    meaning $\|w^{k}\|$ is $\Theta(\|u^{k}\|^{2})$. 

    This implies that $\lim_{k \to \infty} \frac{x^{k}}{\|x^{k}\|} = \tilde{d}$, 
    so by Lemma~\ref{lem:limit} applied to each constraint, $\tilde{d}$ is a recession 
    direction of $S$. We further have that 
    \[q_{i}^{\T} \tilde{d} = \lim_{k \to \infty} \frac{q_{i}^{\T} w^{k}}{\|w^{k}\|} 
    \leq \limsup_{k \to \infty} -\frac{1}{2} a_{i} \frac{\|u^{k}\|^{2}}{L(1 + \|u^{k}\|^{2})} 
    = -\frac{a_{i}}{2L} < 0 .\]
    Since $i \in I$ was arbitrary, we have $q_{i}^{\T} \tilde{d} < 0$ for each $i \in I$.

    That is, for every constraint for which moving in the direction $e$ causes a quadratic increase, 
    moving in the direction $\tilde{d}$ causes those constraints to decrease. Our next step 
    is to replace these directions with rational directions with the same properties. 
    For $\tilde{d}$, since the recession cone of $S$ is given by 
    rational inequalities by Lemma~\ref{lem:convex_recession}, 
    we can apply Lemma~\ref{lem:rational_solution} to get a rational 
    recession direction $d$ of $S$, and by taking $\epsilon$ small 
    enough we still have $q_{i}^{\T} d < 0$ for each $i \in I$.

    Let $J := \{ i \in \{1,2,\dots,m\} : q_{i}^{\T} d = 0 \}$. These are the constraints 
    where moving along $d$ does not change the value, and note 
    that it is disjoint from $I$. We want to replace $e$ with a direction 
    which is both rational, and does not behave badly on the constraints indexed by $J$.
    Define the cone 
    \[\mathcal{C} := \left\{ \rho = \sum_{i \in J} \mu_{i} q_{i} : \mu \geq 0, \pi_{N}(\rho) = 0 \right\} \subseteq U .\]
    Now we have that for each $x \in S$ the inequality $q_{i}^{\T} x \leq -c_{i}$ holds 
    for each $i=1,2,\dots,m$, because the quadratic term is always nonnegative. Then for any $\rho \in \mathcal{C}$ and any $k$,
    we can write 
    \[\rho^{\T} x^{k} = \rho^{\T} u^{k} \leq -\sum_{i \in J} \mu_{i} c_{i}, \]
    since $\rho^{\T} w^{k} = 0$ by definition of $\mathcal{C}$. Dividing this inequality by $\|u^{k}\|$ and taking the limit as $k \to \infty$, we get $\rho^{\T} e \leq 0$, so $e \in \mathcal{C}^{\circ}$, the polar cone of $\mathcal{C}$. 

    This is the property we wanted $e$ to have, and now we show how to find a rational direction close 
    to $e$ which keeps the desired properties.
    $\mathcal{C}$ can be expressed as the image of $\{\mu \geq 0 : \pi_{N}\left(\sum_{i \in J} \mu_{i} q_{i} \right) = 0 \}$ 
    under a rational linear map, so $\mathcal{C}$ is a rational polyhedral cone. This means it has a finite 
    set of rational generators $\rho_{1}, \dots, \rho_{r}$. Then we have
    \[\mathcal{C}^{\circ} = \{x \in \R^{n} : \rho_{i}^{\T} x \leq 0 \quad i=1,2,\dots,r \} ,\]
    which is a rational polyhedron. Also consider $Z = U \cap \bigcap_{i \notin I} \ker Q_{i}$, 
    and we have $e \in Z$, and want to preserve this as well. Then we can 
    apply Lemma~\ref{lem:rational_solution} to $\mathcal{C}^{\circ} \cap Z$ to get a rational 
    element $e' \in \mathcal{C}^{\circ} \cap Z$, and by taking $\epsilon$ small enough 
    we still have $e'^{\T} Q_{i} e' > 0$ for $i \in I$, and the point of intersecting 
    with $Z$ is that we ensure we still have $Q_{i} e' = 0$ for $i \notin I$.

    Now, we want to adjust $e'$ to be well behaved on each constraint indexed by $J$, since 
    those are the constraints where $d$ does not cause a decrease.
    Towards this end, we want to find an adjustment vector $\delta$ such that $q_{i}^{\T} (e' + \delta) \leq 0$ 
    for each $i \in J$. Let $V_{N}$ be a matrix whose columns are a rational basis for $N$, so we can write the 
    condition as finding a vector $\omega$ such that $q_{i}^{\T} V_{N} \omega \leq -q_{i}^{\T} e'$ for each $i \in J$. 
    By Farkas' lemma, this system is infeasible iff there is $\mu \geq 0$ with 
    \[\sum_{i \in J} \mu_{i} V_{N}^{\T} q_{i} = 0, \sum_{i \in J} -\mu_{i} q_{i}^{\T} e' < 0. \]
    But the first condition says that $\rho = \sum_{i \in J} \mu_{i} q_{i}$ has $\pi_{N}(\rho) = 0$, 
    while the second says $\rho^{\T} e' > 0$, which contradicts that $e' \in \mathcal{C}^{\circ}$, 
    so the system is actually feasible, and has a rational solution $\bar{\omega}$, 
    and we can then take $\delta = V_{N} \bar{\omega}$ as our adjustment vector, 
    and define $\bar{e} = e' + \delta$.

    Since $\delta \in N$, we still have $\bar{e}^{\T} Q_{i} \bar{e} > 0$ iff $i \in I$, and 
    now $q_{i}^{\T} \bar{e} \leq 0$ for $i \in J$, while for $d$ we have $q_{i}^{\T} d < 0$ 
    for $i \notin J$, and $I \cap J = \emptyset$.

    Now to actually construct our counterexample using these directions, fix $\bar{x} \in \operatorname{int} S \cap \Q^{n}$, 
    which exists since $S$ is full-dimensional,
    and let $v \in \Q^{n} \setminus \operatorname{span}(\bar{e}, d)$. Fix a rational $\epsilon$ 
    such that $\bar{x} + yv \in S$ for all $|y| \leq \epsilon$. Let 
    \[H := \{x(s,t,y) = \bar{x} + s\bar{e} + td + yv : s,t,y \in \R \}. \]
    Note that $H$ is 3 dimensional: $\bar{e}$ and $d$ are linearly independent, since $\bar{e}^{\T} Q_{i} \bar{e} > 0$ 
    while $Q_{i} d = 0$ for each $i \in I$.
    Now we can express our constraints on $H$ as 
    \begin{align*}
        &i \in I: &f_{i}(x(s,t,y)) = \alpha_{i} s^{2} + \beta_{i}(y) s + f_{i}(\bar{x}+yv) + t q_{i}^{\T} d \\
        &i \notin I: &f_{i}(x(s,t,y)) = f_{i}(\bar{x} + yv) + s q_{i}^{\T} \bar{e} + t q_{i}^{\T} d 
    \end{align*}
    Where $\alpha_{i} = \bar{e}^{\T} Q_{i} \bar{e}$ and $\beta_{i}(y) = 2 (\bar{x} + yv)^{\T} Q_{i} \bar{e} + q_{i}^{\T} \bar{e}$.
    Then if we define for $i \in I$ the function
    \[\psi_{i}(s,y) = -\frac{\alpha_{i} s^{2} + \beta_{i}(y) s + f_{i}(\bar{x} + yv)}{q_{i}^{\T} d} , \] 
    which is well defined since $q_{i}^{\T} d < 0$ for $i \in I$, 
    we have $f_{i}(x(s,t,y)) \leq 0$ iff $t \geq \psi_{i}(s,y)$. Thus, we note 
    \[S \cap H \subseteq \left\{ x(s,t,y) : t \geq \max_{i \in I} \psi_{i}(s,y) \right\} . \]

    Now, we order the functions $\psi_{i}$, $i \in I$. Let $\gamma_{i} := -\frac{1}{q_{i}^{\T} d} > 0$, so that
    \[\psi_{i}(s,y) = \gamma_{i}\alpha_{i} s^{2} + \gamma_{i}\beta_{i}(y) s + \gamma_{i} f_{i}(\bar{x}+yv).\]
    For a fixed value of $y$, consider the triple
    $\gamma_{i}\left( \alpha_{i}, \beta_{i}(y), f_{i}(\bar{x} + yv)\right)$,
    and order these triples lexicographically. For two distinct indices $i,j \in I$, the affine functions
    $\gamma_{i}\beta_{i}(y)$ and $\gamma_{j}\beta_{j}(y)$ are equal for at most one value of $y$, unless they are
    equal everywhere. Similarly, the quadratic functions $\gamma_{i} f_{i}(\bar{x}+yv)$ and $\gamma_{j} f_{j}(\bar{x}+yv)$
    are equal for at most two values of $y$, unless they are equal everywhere. Thus, we can find an interval
    $[y_{1}, y_{2}] \subseteq (-\epsilon, \epsilon)$ with $y_{1} < y_{2}$ that avoids these finitely many values.
    Then the lexicographic order of the triples is the same for all $y \in [y_{1}, y_{2}]$,
    where two indices whose triples are equal for every $y$ are tied.

    Let $i^{*} \in I$ be an index whose triple is largest in this order. We claim that there is $s_{0}$ such that
    $\psi_{i^{*}}(s,y) \geq \psi_{i}(s,y)$ for every $i \in I$, $s \geq s_{0}$, and $y \in [y_{1}, y_{2}]$.
    Indeed, $\psi_{i^{*}}(s,y) - \psi_{i}(s,y) = a s^{2} + b(y) s + c(y)$, where $a$ is a constant,
    $b$ is affine, and $c$ is quadratic. If $a$, $b$, and $c$ are identically zero on $[y_{1}, y_{2}]$, the claim is clear.
    Otherwise, the first of them that is not identically zero is positive on $[y_{1}, y_{2}]$. Since it is continuous
    on the compact interval $[y_{1}, y_{2}]$, it is bounded away from 0, while the later coefficients are bounded.
    This gives $s_{0}$.

    We show that for all sufficiently large $s$, every $y \in [y_{1}, y_{2}]$, and every $t \geq \psi_{i^{*}}(s,y)$,
    we have $x(s,t,y) \in S$. We do this by cases on the constraint.
    \begin{itemize}
        \item $i \in I$, then $f_{i}(x(s,t,y)) = q_{i}^{\T} d (t - \psi_{i}(s,y)) \leq q_{i}^{\T} d (\psi_{i^{*}}(s,y) - \psi_{i}(s,y))$ 
            using that $q_{i}^{\T} d < 0$ for $i \in I$. By the claim, this is $\leq 0$ for $s \geq s_{0}$,
            showing $x(s,t,y)$ satisfies this constraint.
        \item $i \in J$, then $f_{i}(x(s,t,y)) = f_{i}(\bar{x} + yv) + s q_{i}^{\T} \bar{e}$ since $q_{i}^{\T} d = 0$ 
            for $i \in J$. Then $f_{i}(\bar{x}+yv) \leq 0$ since $|y| \leq \epsilon$, and since $q_{i}^{\T} \bar{e} \leq 0$ 
            for $i \in J$, for $s \geq 0$ we have $x(s,t,y)$ feasible for this constraint.
        \item $i \in \{1,2,\dots,m\} \setminus (I \cup J)$. Then
            $f_{i}(x(s,t,y)) = f_{i}(\bar{x} + yv) + s q_{i}^{\T} \bar{e} + t q_{i}^{\T} d$ 
            and $q_{i}^{\T} d < 0$. Since $\gamma_{i^{*}}\alpha_{i^{*}} > 0$ and the other two coefficients of $\psi_{i^{*}}$
            are bounded on $[y_{1},y_{2}]$, there are $\kappa > 0$ and $s_{1}$ such that $\psi_{i^{*}}(s,y) \geq \kappa s^{2}$
            for all $s \geq s_{1}$ and $y \in [y_{1},y_{2}]$. Since $f_{i}(\bar{x}+yv) \leq 0$ and $t \geq \kappa s^{2}$,
            we obtain $f_{i}(x(s,t,y)) \leq s q_{i}^{\T}\bar{e} + \kappa s^{2} q_{i}^{\T} d$,
            which is negative for all sufficiently large $s$.
    \end{itemize}
    Therefore, we have shown that there is some $\bar{s} \in \R$ such that
    \begin{align}
    \label{eq:containment}
    \begin{split}
        \{x(s,t,y) : s \geq \bar{s}, y \in [y_{1}, y_{2}], t \geq \psi_{i^{*}}(s,y) \} 
        \subseteq S \cap H \subseteq \{x(s,t,y) : t \geq \psi_{i^{*}}(s,y) \} .
    \end{split}
    \end{align}
    
    With all of this established, we now construct our objective function and mixed integer lattice. Let $\theta \in [y_{1}, y_{2}]$ be 
    a rational number in the interior of the interval, and $\pi_{H}$ denote the orthogonal projection onto $H$. Then for any $x \in \R^{n}$, 
    we can define $(s(x), t(x), y(x))$ as the unique triple such that $x(s(x),t(x),y(x)) = \pi_{H}(x)$. We 
    can then define our objective as 
    \[f_{0}(x) = t(x) - \psi_{i^{*}}(s(x),y(x)) + (y(x)-\theta)^{2}. \]
    In other words, $f_{0}$ is defined through its restriction to $H$.
    The function $f_{0}$ is quadratic with rational coefficients. Indeed,
    $(s(x),t(x),y(x)) = (V_{H}^{\T}V_{H})^{-1}V_{H}^{\T}(x - \bar{x})$, where $V_{H}$ is the rational matrix with
    columns $\bar{e}$, $d$, and $v$. Moreover, $\psi_{i^{*}}$ has rational coefficients, since $\bar{x}$, $v$, $\bar{e}$, $d$,
    and $f_{i^{*}}$ are rational, and $\theta$ is rational.
    Now $f_{0}(x) \geq 0$ on $S \cap H$, since we showed $t \geq \psi_{i^{*}}(s,y)$ for each $x(s,t,y) \in S \cap H$. 
    We construct a mixed integer lattice so that we can never attain 0, but are still able to approach it. 
    Specifically, we note that the objective is 0 iff $y = \theta$ and $t = \psi_{i^{*}}(s,y)$, so we want to make 
    this impossible.

    Next, we write
    \[\psi_{i^{*}}(s,y) = As^{2} + B_{1}sy + B_{0}s + C(y), \]
    where $A,B_{1},B_{0} \in \Q$, $A = \gamma_{i^{*}}\alpha_{i^{*}} > 0$, and $C$ is a quadratic function with rational
    coefficients. Since $A > 0$, we can choose $\lambda \in \{1,2\}$ such that $2A \lambda + B_{1} \neq 0$.
    Let $D$ be the least common denominator of the coefficients of $\psi_{i^{*}}(k + \lambda \theta, \theta)$ as a
    polynomial in $k$, which are rational since $\lambda$ and $\theta$ are rational. The lattice is given by
    \[\Lambda_{H} = \{x(k + \lambda \sigma, j + \frac{1}{2D}, \sigma) : k,j \in \Z, \sigma \in \R \}. \]
    We extend $\Lambda_{H}$ to a rational mixed integer lattice $\Lambda$ in $\R^{n}$ by adding $n-3$
    linearly independent rational directions orthogonal to the linear subspace parallel to $H$,
    with $p-2$ of them integer and $n-p-1$ of them continuous. Then $\Lambda$ has $p$ integer components and $\Lambda \cap H = \Lambda_{H}$.
    
    We must have $f_{0}(x) > 0$ on $S \cap H \cap \Lambda$, because in order for it 
    to be equal to 0, we must have $\sigma = \theta$ and then $j + \frac{1}{2D} = \psi_{i^{*}}(k+\lambda \theta, \theta)$. 
    But this is impossible, since $\psi_{i^{*}}(k + \lambda \theta, \theta)$ is an integer 
    multiple of $\frac{1}{D}$ since $k \in \Z$, while $j + \frac{1}{2D}$ cannot be an integer 
    multiple of $\frac{1}{D}$.

    To show that we can approach 0, for $k \in \Z$ define $h_{k}(\sigma) = \psi_{i^{*}}(k + \lambda \sigma, \sigma)$,
    and compute the derivative
    \[h_{k}'(\sigma) = (2A \lambda + B_{1}) k + 2A \lambda^{2} \sigma + 2 B_{1} \lambda \sigma
        + \lambda B_{0} + C'(\sigma) . \]
    Let $\nu := |2A \lambda + B_{1}| > 0$. For all $k$ large enough, we have
    $|h_{k}'(\sigma)| \geq \frac{\nu k}{2}$ for all $\sigma \in [y_{1},y_{2}]$, since
    the other terms are bounded independently of $k$.

    Let $\epsilon_{k} = \frac{4}{\nu k}$, so for sufficiently large $k$ we 
    have $[\theta-\epsilon_{k}, \theta+\epsilon_{k}] \subseteq [y_{1},y_{2}]$ 
    and also $k + \lambda \sigma \geq \bar{s}$ for $\sigma \in [\theta-\epsilon_{k}, \theta+\epsilon_{k}]$.
    Now, the image $h_{k}([\theta-\epsilon_{k}, \theta+\epsilon_{k}])$ is an interval
    since $h_{k}$ is continuous. By the mean value theorem and the lower bound on $|h_{k}'|$,
    this interval has length at least $2 \epsilon_{k} \frac{\nu k}{2} = 4$,
    so it contains some $t_{k} \in \Z + \frac{1}{2D}$. Let $\sigma_{k} \in [\theta-\epsilon_{k}, \theta+\epsilon_{k}]$
    be such that $h_{k}(\sigma_{k}) = t_{k}$, and define
    \[z^{k} = x(k + \lambda \sigma_{k}, t_{k}, \sigma_{k}). \]
    Then $z^{k} \in \Lambda$, and also $z^{k} \in S$ for sufficiently 
    large $k$ by~\eqref{eq:containment}. We can compute 
    \begin{align*}
        f_{0}(z^{k}) &= t_{k} - \psi_{i^{*}}(k + \lambda \sigma_{k}, \sigma_{k}) + (\sigma_{k} - \theta)^{2} \\
        &= (\sigma_{k} - \theta)^{2} \leq \epsilon_{k}^{2} , 
    \end{align*}
    which tends to 0 as $k \to \infty$, showing the infimum is indeed 0, but not attained. 
\end{proof}

\begin{proof}[of Corollary~\ref{cor:low_dimension}]
    By definition, every $\text{FW}_{n}$ set is $\text{FW}_{0}$. Conversely, assume that $S$ is $\text{FW}_{0}$.
    Let $\Lambda$ be a rational mixed integer lattice, let $P$ be a rational polyhedron, and let $f_{0}$ be a quadratic
    function with rational coefficients that is bounded from below on $S \cap P \cap \Lambda \neq \emptyset$.
    After the change of variables $x = \tau(y)$, we may assume that $\Lambda = \Z^{p} \times \R^{n-p}$.
    The transformed set is still $\text{FW}_{0}$, since $\tau$ is a rational affine bijection, and its
    affine hull $A$ is still rational, of dimension $k \leq 2$.

    We apply~\cite[Theorem 1.1]{DelPia24} to the polyhedron $A$. Since $A \cap \Lambda$ is nonempty,
    it gives $p' \leq p$, $n'$, and a map $\sigma(z) = \tilde{x} + Gz$ with $\tilde{x} \in \Q^{n}$ and
    $G \in \Q^{n \times n'}$ of full column rank, such that $A = \sigma(A')$ for a full-dimensional
    polyhedron $A' \subseteq \R^{n'}$, and $A \cap \Lambda = \sigma(A' \cap \Lambda')$, where
    $\Lambda' := \Z^{p'} \times \R^{n'-p'}$. Since $A$ is an affine subspace and $\sigma$ is injective,
    we have $A' = \R^{n'}$ and $n' = k$.

    Let $S' := \{z \in \R^{k} : f_{i}(\sigma(z)) \leq 0, \ i = 1,2,\dots,m\}$, so that $S = \sigma(S')$.
    Let $G^{+} := (G^{\T} G)^{-1} G^{\T}$. The set $S'$ is $\text{FW}_{0}$. Indeed, for a rational
    polyhedron $P'' \subseteq \R^{k}$ and a quadratic function $g$ with rational coefficients, the map
    $\sigma$ is a bijection from $S' \cap P''$ to $S \cap \sigma(P'')$, the set
    $\sigma(P'') = \{x \in A : G^{+}(x - \tilde{x}) \in P''\}$ is a rational polyhedron, and
    $g(z) = g(G^{+}(\sigma(z) - \tilde{x}))$. Since $k \leq 2$, Corollary~\ref{cor:all_integer} implies that
    $S'$ is $\text{FW}_{k}$, where for $k = 0$ there is nothing to prove.
    Finally, $S \cap P \cap \Lambda = \sigma(S' \cap P' \cap \Lambda')$, where $P' := \sigma^{-1}(P)$ is a
    rational polyhedron, and $f_{0} \circ \sigma$ is a quadratic function with rational coefficients.
    Hence $f_{0}$ attains its infimum on $S \cap P \cap \Lambda$.
\end{proof}

\section{Conditions on the objective}
\label{app:convex}

This appendix contains the proof of Theorem~\ref{thm:convex} and Corollary~\ref{cor:one_integer_quasi-convex}.

The following lemma is the inductive step of our argument: it reduces the dimension by one.
When we project out a recession direction $d$, the inequalities that do not depend on $d$ are kept.
The other inequalities can be dropped, because adding a large multiple of $d$ to any point satisfies them.

\begin{lemma} \label{lem:lower_dimension}
    Consider an instance of~\eqref{prob:MIQCQP} with rational data, where
    all constraints $f_{i}$ $(i=1,2,\dots,m)$ are convex and $f_{0}$ is bounded 
    from above on $S \cap \Lambda$.
    Let $d \in \Q^{n}$ be a nonzero recession direction for all constraints, 
    such that $f_{0}(x) = f_{0}(x+td)$ for every $x \in S$ 
    and $t \in \R$.

    Then there exist $E \in \Q^{n \times (n-1)}$, $m' \in \Q^{n}$, $p' \in \{p-1, p\}$, and a rational mixed integer
    lattice $\Lambda' \subseteq \R^{n-1}$ with $p'$ integer components such that the problem
    \begin{align}
        \label{prob:lower_dimension}
        \begin{split}
            \inf\; & f_{0}(E x + m')   \\
            \text{s.t.}\; & x \in S' \cap \Lambda',
        \end{split}
    \end{align}
    where $I := \{k \in \{1,2,\dots,m\} : q_{k}^{\T} d = 0 \}$ and
    \[S' := \{x \in \R^{n-1} : f_{k}(E x + m') \leq 0, \quad k \in I \},\]
    has the same infimum as the original instance of~\eqref{prob:MIQCQP}, and $f_{0}(Ex + m')$ is bounded
    from above on $S' \cap \Lambda'$.
    Moreover, the functions $f_{k}(E x + m')$, $k \in I$, are convex.
    Finally, for every $x \in S' \cap \Lambda'$ there is $\bar{x} \in S \cap \Lambda$ with $f_{0}(\bar{x}) = f_{0}(Ex + m')$.
\end{lemma}

\begin{proof} 
    Without loss of generality, we can assume $\Lambda = \Z^{p} \times \R^{n-p}$,
    because we can replace $x$ with $\tau(y)$ and all assumptions are still met.
    Below, we also apply a rational linear change of variables that maps $\Z^{p} \times \R^{n-p}$ onto itself.
    In the new variables, we construct the required map with $m' = 0$. Composing it with the two changes
    of variables gives the rational affine map $x \mapsto Ex + m'$ for the original instance.

    Next, we bring $d$ into a convenient form.
    Let $d = (d_{\Z}, d_{\R})$, scaled by a positive rational number so that $d_{\Z} \in \Z^{p}$ is primitive 
    and $d_{\R} \in \Q^{n-p}$. 
    If $d_{\Z} \neq 0$, let $U \in \Z^{p \times p}$ be a unimodular matrix such that $Ue^{1} = d_{\Z}$. Then 
    consider $T = \begin{pmatrix} U & 0 \\ D & I_{n-p} \end{pmatrix}$ 
    where $D$ has $d_{\R}$ as the first column and 0 for the rest of the entries, 
    and $I_{n-p}$ is the identity matrix. If $d_{\Z} = 0$, then $d_{\R}$ is nonzero, 
    and instead let $U \in \Q^{(n-p) \times (n-p)}$ be an invertible matrix such that $Ue^{1} = d_{\R}$. 
    Then we can take $T = \begin{pmatrix} I_{p} & 0 \\ 0 & U \end{pmatrix}$, 
    since no off-diagonal block is needed when $d_{\Z} = 0$.

    In either case, $T$ is a bijection from $\Z^{p} \times \R^{n-p}$ to itself,
    so we can do a change of basis $x=Ty$ such that the 
    lattice is still $\Z^{p} \times \R^{n-p}$, but our recession direction 
    is now either $e^{1}$ or $e^{p+1}$, and all assumptions are still met.

    Let $E \in \Q^{n \times (n-1)}$ be the matrix that inserts a zero component in position $1$
    or $p+1$, respectively, and let $\Lambda' := \Z^{p'} \times \R^{n-1-p'}$, where $p' = p-1$
    if we remove an integer component and $p' = p$ if we remove a continuous component.

    Now we prove that~\eqref{prob:MIQCQP} and~\eqref{prob:lower_dimension} 
    have the same infimum. 
    Let $x \in S' \cap \Lambda'$, which has objective value $f_{0}(E x)$ for \eqref{prob:lower_dimension}. 
    Consider the corresponding point $Ex \in \R^{n}$.

    By Lemma~\ref{lem:convex_recession}, we have $Q_{i} d = 0$ and $q_{i}^{\T} d \leq 0$ for each $i=1,2,\dots,m$.
    We already know that for every $k \in I, f_{k}(Ex) \leq 0$ since $x \in S' \cap \Lambda'$, 
    while for every $k \in \{1,2,\dots,m\} \setminus I, Q_{k} d = 0, q_{k}^{\T} d < 0$.
    Thus,
    \[\forall k \in \{1,2,\dots,m\} \setminus I, \exists t_{0} > 0, \forall t \geq t_{0}, f_{k}(Ex + td) \leq 0.\]
    Since $d$ is the unit vector of the removed component, $Ex + td \in \Z^{p} \times \R^{n-p}$
    for every $t \in \Z$. Moreover, $f_{k}(Ex + td) = f_{k}(Ex) \leq 0$ for every $k \in I$,
    since $Q_{k} d = 0$ and $q_{k}^{\T} d = 0$.
    Hence $Ex + td \in S \cap \Lambda$ for all sufficiently large $t \in \Z$.
    Since $Ex + td \in S$, the assumption implies that $f_{0}$ is constant on the line through
    $Ex + td$ with direction $d$. This line contains $Ex$, so $f_{0}(Ex + td) = f_{0}(Ex)$.
    Hence, for every such $t$, the point $\bar{x} := Ex + td$ lies in $S \cap \Lambda$ and satisfies
    $f_{0}(\bar{x}) = f_{0}(Ex)$. In particular, $f_{0}(Ex)$ is bounded from above
    on $S' \cap \Lambda'$, since $f_{0}$ is bounded from above on $S \cap \Lambda$.

    For the other direction, take $x \in S \cap \Lambda$ and consider projecting out the chosen 
    component. We can write the resulting point as  
    $x'$ for some $x' \in \Z^{p'} \times \R^{n-1-p'}$. 
    Then $x' \in S' \cap \Lambda'$, because $Ex' = x + td$ for some $t \in \R$, and 
    since $d$ is a recession direction, Lemma~\ref{lem:convex_recession} 
    implies $f_{k}(x) = f_{k}(x+td)$ for every $k \in I$, so 
    $f_{k}(x) \leq 0$ iff $f_{k}(Ex') \leq 0$, 
    and $f_{0}(Ex') = f_{0}(Ex' + td)$ by assumption, proving the infimum is the same.

    Finally, convexity is preserved as $f_{i}(Ex)$ is the composition of a convex and linear 
    function.
\end{proof}

The next two lemmas provide the remaining ingredients for Theorem~\ref{thm:convex}.

\begin{lemma} \label{lem:rational_recession}
    Consider an instance of~\eqref{prob:MIQCQP} with rational data, where all constraints
    $f_{i}$ $(i=1,2,\dots,m)$ are convex. If $S \cap \Lambda$ is unbounded, then there is some $d \in \Q^{n} \setminus \{0\}$ 
    which is a recession direction of all constraints.
\end{lemma}

\begin{proof}
    Let $(x^{k})$ be a sequence of points in $S \cap \Lambda$ with $\lim_{k \to \infty} \|x^{k}\| = \infty$. 
    Replace $(x^{k})$ with a subsequence such 
    that $\left( \frac{x^{k}}{\|x^{k}\|} \right)$ 
    has a well defined limit, which we can do by compactness. 
    Denote $d' = \lim_{k \to \infty} \frac{x^{k}}{\|x^{k}\|}$, 
    and by applying Lemma~\ref{lem:limit} to each constraint, 
    $d'$ must be a recession direction of each of them. 

    By Lemma~\ref{lem:convex_recession}, the set of recession directions 
    of all constraints is given by a polyhedron, which is rational since 
    we assume all functions have rational coefficients. Since $\|d'\|=1$, 
    we can apply Lemma~\ref{lem:rational_solution} to find 
    some $d \in \Q^{n} \setminus \{0\}$ which is a recession direction 
    of all constraints.
\end{proof}

\begin{lemma} \label{lem:added_constraints} 
    Consider an instance of~\eqref{prob:MIQCQP} with rational data, where all constraints
    $f_{i}$ $(i=1,2,\dots,m)$ are convex. Let $d \in \Q^{n} \setminus \{0\}$ be a recession 
    direction of all constraints such that 
    $f_{0}(x) = f_{0}(x+td)$ for every $x \in S \cap \Lambda$ and $t \in \R$. 

    Then if we add $\pm \left( 2 d^{\T} Q_{0} x + q_{0}^{\T} d \right) \leq 0$ as constraints,
    $S \cap \Lambda$ remains unchanged, but now $f_{0}(x) = f_{0}(x+td)$ 
    for every $x \in \tilde{S}$ and $t \in \R$, 
    where $\tilde{S}$ is the continuous relaxation of the feasible set given by 
    \[\tilde{S} := \{x \in \R^{n} : 2d^{\T} Q_{0} x + q_{0}^{\T} d = 0, f_{i}(x) \leq 0, \quad i=1,2,\dots,m \}. \]
    Additionally, $d$ is still a recession direction of these constraints, and they are both convex and rational.
\end{lemma}

\begin{proof}
    By expansion we have 
    \[f_{0}(x+td) = t^{2} d^{\T} Q_{0} d + 2t d^{\T} Q_{0} x + t q_{0}^{\T} d + x^{\T} Q_{0} x + q_{0}^{\T} x + c_{0}. \]
    This is a constant function 
    of $t$ iff $d^{\T} Q_{0} d = 0$ and $2d^{\T} Q_{0}x + q_{0}^{\T} d = 0$. 

    By assumption $f_{0}(x+td)$ is constant as a function of $t$ for each $x \in S \cap \Lambda$, 
    so we must have $d^{\T} Q_{0} d = 0$ and for every $x \in S \cap \Lambda, 2 d^{\T} Q_{0} x + q_{0}^{\T} d = 0$, 
    so adding the additional constraints does not change $S \cap \Lambda$. On the other hand, adding the constraint 
    ensures that these two conditions are met for every $x \in \tilde{S}$, so $f_{0}(x+td)$ is constant 
    as a function of $t$ for each $x \in \tilde{S}$.

    The direction $d$ remains a recession direction of the new constraints, since $d^{\T} Q_{0} d = 0$.
    The new constraints are convex, since they are affine, and they have rational coefficients,
    since $d$, $Q_{0}$, and $q_{0}$ are rational.
\end{proof}

\begin{proof}[of Theorem~\ref{thm:convex}]
    We proceed by induction on the dimension of 
    the problem. If $n=0$, then there is only a single 
    feasible point since we assume the feasible region is nonempty,
    and it trivially attains the infimum of 
    all feasible points.
    
    For $n > 0$, 
    if $S \cap \Lambda$ is bounded, since it is closed as the intersection 
    of the closed sets $S$ and $\Lambda$, it is compact, so by the 
    extreme value theorem $f_{0}$ attains a minimum value on it.

    Otherwise $S \cap \Lambda$ is unbounded, and we apply Lemma~\ref{lem:rational_recession} 
    to get a rational recession direction $d$ of all constraints. 
    By assumption, $f_{0}$ is bounded both above and below on $S \cap \Lambda$. Therefore, 
    for every $x \in S \cap \Lambda, f_{0}(x+td)$ must be a constant function of $t$, because 
    $f_{0}(x+td)$ is a polynomial in $t$. Lemma~\ref{lem:mixed_integer_spacing} 
    implies there are arbitrarily large values of $t$ such that $x+td \in S \cap \Lambda$, 
    so if $f_{0}(x+td)$ were not a constant function of $t$, it would tend 
    to either $\infty$ or $-\infty$ for large values of $t$, contradicting that it 
    is bounded in both directions.

    Then if we add the constraints
    $\pm \left( 2 d^{\T} Q_{0} x + q_{0}^{\T} d \right) \leq 0$ to our system, 
    Lemma~\ref{lem:added_constraints} ensures we do not change $S \cap \Lambda$, 
    but now we meet the conditions 
    to apply Lemma~\ref{lem:lower_dimension}.
    This gives a lower dimensional problem with the same infimum, which by the induction 
    hypothesis attains its infimum at some feasible 
    point $x \in S' \cap \Lambda'$, and as noted in 
    Lemma~\ref{lem:lower_dimension}, there is a corresponding 
    point $\bar{x} \in S \cap \Lambda$ which 
    attains the infimum of~\eqref{prob:MIQCQP}.
\end{proof}

\begin{remark} \label{rmk:convex_real}
    If $p = 0$, that is, $\Lambda = \R^{n}$, then Theorem~\ref{thm:convex} holds without the
    rationality assumptions. Indeed, rationality is used in Lemmas~\ref{lem:mixed_integer_spacing},
    \ref{lem:lower_dimension}, and~\ref{lem:rational_recession} to stay in the mixed integer lattice,
    and in Lemma~\ref{lem:added_constraints} to keep the added constraints rational for the next step
    of the induction. For $p = 0$, the limit direction in the proof of Lemma~\ref{lem:rational_recession}
    can be used directly as the recession direction $d$, the change of variables in the proof of
    Lemma~\ref{lem:lower_dimension} can be real, the points $x + td$ with $t \geq 0$ are in $\Lambda$,
    and the added constraints need not be rational.
\end{remark}

For at most one integer variable, the continuous result used in Theorem~\ref{thm:generic_one_integer}
is Theorem~\ref{thm:convex} with $p = 0$. This gives the following corollary.
\begin{corollary} \label{cor:one_integer_quasi-convex}
    Consider a problem of the form~\eqref{prob:generic_one_integer}, or its continuous version with
    $\Lambda = \R^{n}$, where all $g_{j}$ $(j=0,1,\dots,m)$ are quadratic functions, all constraint functions
    $g_{i}$ $(i=1,2,\dots,m)$ are convex, and $g_{0}$ is quasi-convex over
    $K := \{x \in \R^{n} : g_{i}(x) \leq 0, \ i = 1,2,\dots,m\}$.
    Then the problem attains its infimum (when it is finite).
\end{corollary}

\begin{proof}
    We first show that, for every set $K' \subseteq K$ defined by convex quadratic inequalities,
    the problem $\inf \{g_{0}(x) : x \in K'\}$ attains its infimum when it is finite.
    For $K' = K$, this is the continuous case.
    Let $\bar{y} + W z$, $z \in \R^{k}$, be a parametrization of the affine hull of $K'$,
    where $\bar{y} \in \R^{n}$ and $W \in \R^{n \times k}$ has full column rank.
    In the variables $z$, the feasible set is full-dimensional and defined by convex quadratic
    inequalities. The objective $\tilde{g}_{0}(z) := g_{0}(\bar{y} + Wz)$ is a quadratic function that is
    quasi-convex over the feasible set, since $K' \subseteq K$. Since all variables are continuous,
    this change of variables does not affect attainment. Let $\tilde{Q}$ be the symmetric matrix
    of the quadratic part of $\tilde{g}_{0}$.
    If $\tilde{Q}$ is positive semidefinite, we add the constraint $\tilde{g}_{0}(z) - r \leq 0$
    for some $r$ larger than the infimum, as in the proof of Corollary~\ref{cor:quasi-convex}.
    Otherwise, \cite[Theorem 6.3]{GeneralizedConcavity10} shows that $\tilde{g}_{0}$ is bounded from above
    on the feasible set, as in the proof of Corollary~\ref{cor:quasi-convex}. In both cases, the objective is bounded
    from above on the feasible region, so the infimum is attained by Theorem~\ref{thm:convex}
    with $p = 0$ and Remark~\ref{rmk:convex_real}.

    For the mixed integer case, we apply Theorem~\ref{thm:generic_one_integer}. As in the proof of
    Corollary~\ref{cor:quadratic_one_integer}, the set in its hypothesis is semialgebraic,
    so it has finitely many connected components. The feasible set of every
    relaxation~\eqref{prob:relaxed_with_linear} is $K \cap \{x \in \R^{n} : L \leq a^{\T} x \leq U\}$,
    which is defined by convex quadratic inequalities. Hence every relaxation with a finite infimum attains it.
\end{proof}

\end{document}